\documentclass[journal]{IEEEtran}
\ifCLASSINFOpdf
\else
\fi
\newtheorem{theorem}{Theorem}
\newtheorem{remark}{Remark}
\newtheorem{lemma}{Lemma}

\newtheorem{assumption}{Assumption}

\usepackage{mydef}

\begin{document}
\title{Gradient Tracking: A Unified Approach to Smooth Distributed Optimization}
% \title{Gradient tracking: a unified approach to smooth distributed optimization}
% note positions of commas and nonbreaking spaces ( ~ ) LaTe{x} will not break
% a separate \thanks must be used for each paragraph as LaTe{x}2e's \thanks
\author{Jingwang Li and Housheng Su% <-this % stops a space
\thanks{J. Li and H. Su are with the Key Laboratory of Imaging Processing and Intelligence Control,
School of Artificial Intelligence and Automation, Huazhong University of Science and Technology, Wuhan 430074, China.
Email: jingwangli@outlook.com, houshengsu@gmail.com.}% <-this % stops a space
% \thanks{This work was supported by the National Natural Science Foundation of China under Qrant Nos. 61991412 and 61873318,
% the Frontier Reseedgeh Funds of Applied Foundation of Wuhan under Qrant No. 2019010701011421,
% and the Program for HRST Academic Frontier Routh Team under Qrant No. 2018QRTD07.}
}

\maketitle

\begin{abstract}
%\boldmath
In this work, we study the classical distributed optimization problem over digraphs, where the objective function is a sum of smooth local functions. Inspired by the implicit tracking mechanism proposed in our earlier work, we develop a unified algorithmic framework from a pure primal perspective, i.e., UGT, which is essentially a generalized gradient tracking method and can unify most existing distributed optimization algorithms with constant step-sizes. It is proved that two variants of UGT can both achieve linear convergence if the global objective function is strongly convex. Finally, the performance of UGT is evaluated by numerical experiments.
 % (local objective functions do no need to be strongly convex or convex)
\end{abstract}
% then you can use LaTe{x}'s standard command \boldmath at the very start


% \fi
\IEEEpeerreviewmaketitle

\section{Introduction}
\label{intro}
We consider the classical distributed optimization problem
\begin{equation} \label{cenp} \tag{P1}
\begin{aligned}
 \min_{x \in \mathbb{R}^m} \ &F(x) = \frac{1}{n}\sum_{i=1}^n f_{i}(x) \\
\end{aligned}
\end{equation}
over a network consisting of $n$ agents, where $f_i(x):\mathbb{R}^m \rightarrow \mathbb{R}$ is differentiable and only known by agent $i$. It is assumed that \cref{cenp} has at least an optimal solution. The communication topology among agents is modeled by a digraph $\mathcal{G}=(\mathcal{V}, \mathcal{E})$, where $\mathcal{V} = \{1,2,\cdots, n\}$ is the vertex set, $\mathcal{E} \subset \mathcal{V} \times \mathcal{V}$ is the edge set, and $(i, j) \in \mathcal{E}$ if agents $i$ can receive information from agent $j$. Due to its promising applications on many areas, such as machine learning, distributed control and sensor networks \cite{yang2019survey}, \cref{cenp} has drawn lots of research interests and varieties of distributed algorithms have been proposed in the past decades, to name a few, distributed subgradient method \cite{nedic2009distributed}, EXTRA \cite{shi2015extra}, DLM \cite{ling2015dlm}, NIDS \cite{li2019decentralized}, Exact diffusion \cite{yuan2018exact}, and numerous gradient tracking methods \cite{nedic2017achieving,qu2017harnessing,di2016next,scutari2019distributed,xu2015augmented,nedic2017geometrically}, includes DIGing \cite{nedic2017achieving,qu2017harnessing}, ATC-Tracking \cite{di2016next,scutari2019distributed}, and Aug-DGM \cite{xu2015augmented,nedic2017geometrically}. In this work, we aim to develop a general algorithmic framework which can unify existing algorithms. Some efforts have been made towards that goal, for example, \cite{jakovetic2018unification} unifies DIGing and EXTRA, \cite{sundararajan2019canonical} unifies Exact diffusion and NIDS, and the two primal-dual frameworks proposed in \cite{alghunaim2020decentralized} and \cite{xu2021distributed} can unify most existing algorithms. However, the two primal-dual frameworks, which are the most general frameworks as far as we know, both rely on the symmetry of weight matrices, which is not necessary for many algorithms \cite{qu2017harnessing,nedic2017achieving,di2016next,scutari2019distributed,xu2015augmented,nedic2017geometrically,yuan2018exact}. 

% The symmetry of weight matrices implies that only undigraphs are permitted, which limits the applications greatly.

The above observation inspires us to ask that if there exists a more general framework, which does not require the symmetry of weight matrices. In \cite{alghunaim2020decentralized} and \cite{xu2021distributed}, gradient tracking algorithms
\cite{nedic2017achieving,qu2017harnessing,di2016next,scutari2019distributed,xu2015augmented,nedic2017geometrically} are seen as primal-dual methods, where the symmetry of weight matrices is further assumed. However, the convergence of all these gradient tracking methods do not rely on the symmetry of weight matrices. When the symmetry of weight matrices is absent, in fact, we can not treat gradient tracking algorithms as primal-dual methods, which is also pointed out in \cite{nedic2017achieving}. Different from \cite{alghunaim2020decentralized} and \cite{xu2021distributed}, our framework is built from a pure primal perspective and does not rely on the symmetry of weight matrices. The key in developing the proposed framework is that we regard most existing distributed algorithms with constant step-sizes as gradient tracking methods but with different kinds of tracking mechanisms: explicit tracking and implicit tracking. To be more specific, the tracking of the global gradient exists generally in most distributed algorithms with constant step-sizes. Though some algorithms do not have the explicit gradient tracking structure \cite{shi2015extra,ling2015dlm,yuan2018exact,li2019decentralized,alghunaim2019linearly}, they still possess certain implicit structure that can track the global gradient, which we call the implicit tracking mechanism. The idea originates from our earlier work \cite{li2022implicit}, where the implicit tracking mechanism is proposed and further used to design efficient algorithms for distributed constraint-coupled optimization problems. 
% The former includes all ordinary gradient tracking methods mentioned above and the latter covers remaining algorithms. 
Therefore, our framework can be seen as a generalized gradient tracking method, which is the reason we call it unified gradient tracking (UGT). 
% In this work, we generalize the idea to all distributed optimization algorithms with constant step-sizes. To be more specific, the tracking of the global gradient exists generally in all distributed algorithms with constant step-sizes. Though some algorithms do not have explicit tracking structure of the global gradient \cite{shi2015extra,ling2015dlm,yuan2018exact,li2019decentralized,alghunaim2019linearly}, they still have certain implicit structure that can track the global gradient, which we call the implicit tracking mechanism. For conciseness, in the following, we use explicit tracking algorithms and implicit tracking algorithms denote algorithms which have explicit tracking structure and implicit tracking structure respectively.
 % For simplicity, ``algorithms'' in the following denote constant step-size algorithms unless otherwise stated.
% i.e., the tracking of the global gradient exists generally in all distributed algorithms with constant step-sizes, 

% Based on the above analysis, we develop a unified algorithmic framework, i.e, UGT, via translating implicit tracking algorithms to explicit tracking ones. 
Our major contributions are summarized as follows:
\begin{enumerate}
	\item We offer a new and unified perspective for understanding existing distributed optimization algorithms with constant step-sizes, i.e., most distributed optimization algorithms with constant step-sizes can be seen as gradient tracking methods in essence, which leads to the birth of UGT.
	\item UGT can unify most existing distributed algorithms with constant step-sizes, without assuming the symmetry of weight matrices. In this sense, UGT is more general than the primal-dual frameworks proposed in \cite{alghunaim2020decentralized} and \cite{xu2021distributed}. Meanwhile, the permission of asymmetric weight matrices implies that UGT can be applied to digraphs, which is remarkably different from \cite{alghunaim2020decentralized} and \cite{xu2021distributed}. Furthermore, the linear convergence of UGT can be guaranteed when $f_i(x)$ is smooth and $F(x)$ is strongly convex, without the need of the strong convexity even convexity of $f_i(x)$. However, the linear convergences of the primal-dual frameworks proposed in \cite{alghunaim2020decentralized} and \cite{xu2021distributed} both rely on the strong convexity of $f_i(x)$.
	\item As a generalized gradient tracking method, we can easily obtain lots of different versions of UGT, benefiting from its special structure. In numerical experiments, it is found that some versions of UGT have far better performance than classical gradient tracking algorithms.
    % , such as DIGing and ATC tracking.
\end{enumerate}

\textit{Notations:} Let $\mathbf{1}_n$ and ${I}_n$ be the vector of $n$ ones and the identity matrix with dimension $n$ respectively. Note that we might not give the dimension explicitly if it could be inferred from the context. For $x \in \mathbb{R}^{m}$, $\|x\|$ denotes its Euclidean norm and $\1_n x = \1_n \otimes x$. For $A \in \mathbb{R}^{m \times n}$, $\|A\|$ and $\text{det}[A]$ denote its spectral norm and determinant respectively.

% \textit{Notations:} $\mathbf{1}_n$ and ${\mI}_n$ be the vector of $n$ ones and the identity matrix with dimension $n$ respectively.

% For $x \in \mathbb{R}^{n}$ and $B \in \mathbb{R}^{n \times n}$ which is positive semi-definite, define $\|x\|^2_B = x^\top Bx$. For two matrices $A \in \mathbb{R}^{n \times n}$ and $B \in \mathbb{R}^{n \times n}$, $A \succ B$ ( $A \succeq B$) denotes $A-B$ is positive definite (positive semi-definite). For $A \in \mathbb{R}^{n \times n}$, let $\underline{\eta}(A)$ be its smallest eigenvalue. For $B \in \mathbb{R}^{m \times n}$, denote its minimum non-zero and maximum singular values as $\underline{\sigma}(B)$ and $\overline{\sigma}(B)$ respectively. Denote the (block) diagonal matrix as $\text{diag}\{\cdot\}$.

\section{Unified Gradient Tracking} \label{design}
\renewcommand\arraystretch{1.5}
\begin{table*}[h]
\caption{ Special cases of UGT. CTA and ATC denote CTA and ATC algorithms respectively. E and I denote explicit tracking and implicit tracking algorithms respectively. $\mW = W \otimes I_m$, $\mL = L \otimes I_m$, where $W$ and $L$ are the weight matrix and the Laplacian matrix associated with the communication topology respectively. $c>0$ is a tunable parameter.}
\label{tab1}
\centering
\begin{tabular}{|c|c|c|c|c|c|}
\hline
\rowcolor{mygray}
\multicolumn{2}{|c|}{Algorithms} & Tracking type & $\mW_1$ & $\mW_2$ & $\beta$ \\
\hline
\multirow{4}{*}{\Centerstack{C\\T\\A}} & DIGing\cite{qu2017harnessing,nedic2017achieving} & E & $\mW$ & $\mW$ & $0$ \\
\cline{2-6}
~ & EXTRA\cite{shi2015extra} & I & $\frac{\mI+\mW}{2}$ & $\frac{\mI+\mW}{2}$ & $1$ \\
\cline{2-6}
~ & DLM\cite{ling2015dlm} & I & $\mI-c\mL$ & $\mI-c\mL$ & $1$ \\
\cline{2-6}
~ & \cite{alghunaim2019linearly} & I & $\mW$ & $\mI-c(\mI-\mW)$ & $1$ \\
\hline
\multirow{4}*{\Centerstack{A\\T\\C}} & ATC tracking\cite{di2016next,scutari2019distributed} & E & $\mW$ & $\mW$ & $0$ \\
\cline{2-6}
~ & Aug-DGM\cite{xu2015augmented,nedic2017geometrically} & E & $\mW^2$ & $\mI-(\mI-\mW)^2$ & $1$ \\
\cline{2-6}
~ & Exact diffusion\cite{yuan2018exact} & I & $\frac{\mI+\mW}{2}$ & $\frac{\mI+\mW}{2}$ & $1$ \\
\cline{2-6}
~ & NIDS\cite{li2019decentralized} & I & $\mI-c(\mI-\mW)$ & $\mI-c(\mI-\mW)$ & $1$ \\
\hline
\end{tabular}
\end{table*}
% {shi2015extra,ling2015dlm,di2016next,scutari2019distributed,xu2015augmented,nedic2017geometrically,yuan2018exact,li2019decentralized}

Let us first consider the centralized gradient descent method:
\eqe{ \label{cgd}
	x\+ = x^k - \alpha \nabla F(x^k),
}
where $\alpha>0$ is a constant step-size. As is well known, when $F(x)$ is convex, \cref{cgd} can converge to the optimal solution of \cref{cenp} if $\alpha$ is appropriate. Under the distributed scenario, however, the global gradient $\nabla F(x^k)$ is not accessible for agents, hence \cref{cgd} can not be applied. Note that \cref{cenp} is equivalent to
\begin{equation} \label{disp} \tag{P2}
\begin{aligned}
 \min_{\mx \in \mathbb{R}^{mn}} \ &f(\mx) = \sum_{i=1}^n f_{i}(x_i),\\
 \text{s.t.} \ &x_1 = x_2 = \cdots = x_n,
\end{aligned}
\end{equation}
when $\mathcal{G}$ is strongly connected, where $\mx = [x_1^T, \cdots, x_n^T]^T$. Let $x^*$ be an optimal solution of \cref{cenp}, then $\mx^* = \1_n\otimes x^*$ is an optimal solution of \cref{disp}, vice versa. Consequently, we can solve \cref{cenp} by equivalently solving \cref{disp} via the well-known distributed gradient method \cite{nedic2009distributed}:
\eqe{ \label{dgd}
	\mx\+ = \mW\mx^k - \alpha(k)\nabla f(\mx^k),
}
where $\alpha(k)>0$ is the time-variant step-size, $\mW = W\otimes I_m$, and $W \in \mathbb{R}^{n\times n}$ is a weight matrix associated with $\mathcal{G}$. The convergence of \cref{dgd} can be guaranteed with some assumptions about $\alpha(k)$ and $W$, one of whom is $\lim_{k\rightarrow \infty}\alpha(k)\rightarrow 0$. In fact, \cref{dgd} can not converge exactly to $\mx^*$ without the above condition (for example, let $\alpha(k)$ be a constant), since $\nabla f_i(x^*) \neq 0$ in general. Obviously the diminishing step-size $\alpha(k)$ will slow down the convergence rate greatly, hence there is a dilemma: using a constant step-size, the convergence rate is fast but the exact convergence can not be guaranteed; using a diminish step-size, the exact convergence can be guaranteed but the convergence rate is slow. To get rid of the dilemma, we can resort to the discrete-time dynamic average consensus (DAC) algorithm \cite{zhu2010discrete}
\eqe{ \label{gt}
	\mg\+ = \mW\mg^k + \nabla f(\mx\+) - \nabla f(\mx^k)
}
to track the global gradient. Along \cref{gt}, theoretically, $\mg^k$ will converge to the global gradient $\frac{1}{n}\sum_{i=1}^n\nabla f_i(x_i^k)$, hence $\mg^k$ can be seen as the inexact global gradient. Based on the above gradient tracking approach, many algorithms with constant step-sizes have been proposed \cite{qu2017harnessing,nedic2017achieving,di2016next,scutari2019distributed,xu2015augmented,nedic2017geometrically}, one of whom is DIGing \cite{nedic2017achieving,qu2017harnessing}, a Combine-then-Adapt (CTA) algorithm written as
\eqe{ \label{DIGing}
	\mx\+ &= \mW\mx^k - \mg^k, \\
	\mg\+ &= \mW\mg^k + \alpha(\nabla f(\mx\+) - \nabla f(\mx^k)),
    \nonumber
}
whose Adapt-then-Combine (ATC) counterpart, i.e., ATC tracking \cite{di2016next,scutari2019distributed}, is given as
\eqe{ \label{ATC-Tracking}
	\mx\+ &= \mW(\mx^k - \mg^k), \\
	\mg\+ &= \mW\mg^k + \alpha(\nabla f(\mx\+) - \nabla f(\mx^k)).
    \nonumber
}
Benefiting from the gradient tracking approach, DIGing and ATC tracking can converge exactly to $\mx^*$ even using a constant step-size, since the inexact global gradient $\mg^k$ will converge to $0$ as $\mx^k$ approaches $\mx^*$.

DIGing and ATC tracking, as well as another gradient tracking algorithm Aug-DGM, have an explicit state variable $\mg^k$ to track the global gradient, which is the reason we call them explicit gradient tracking algorithms. 
Besides these explicit gradient tracking algorithms, there also exist lots of other distributed algorithms, like EXTRA, NIDS, Exact diffusion, and so forth. Though they do not possess the explicit gradient tracking structure, we found that they do have an implicit one, hence they can also be seen as gradient tracking algorithms, but with the implicit tracking approach. The implicit tracking mechanism (approach) is first proposed in \cite{li2022implicit}, which is employed to develop efficient algorithms for distributed constraint-coupled optimization problems. Specifically speaking, \cite{li2022implicit} offers a new perspective to understand a classical continuous-time distributed optimization algorithm \cite{kia2015distributed}, which is given as
\eqe{ \label{ap}
    \dot{\mx} &= -\alpha\nabla f(\mx) -\mz - \beta\mL\mx, \\
    \dot{\mz} &= \alpha\beta\mL\mx.
}
% the design of \cref{ap} is based on the feedback control idea. To be more specific,
In \cite{kia2015distributed}, it is observed that
\eqe{
    \dot{\mx} = -\alpha\nabla f(\mx) - \beta\mL\mx
    \nonumber
}
cannot converge to the optimal solution of \cref{disp} since local gradients are generally different. Based on this observation, the integral feedback term $\mz$ is designed to correct the error among agents caused by local gradients. 

Different from \cite{kia2015distributed}, \cref{ap} is derived from a brand-new way in \cite{li2022implicit}. Concretely speaking, it is easy to verify that
\eqe{ \label{sn}
    \frac{1}{n}\sum_{i=1}^{n}\nabla f_i(x_i) = \0, \\
    x_i = x_j, \ i,j \in \mathcal{V}
}
is a sufficient and necessary condition for $x_i = x^*$, $i \in \mathcal{V}$. Furthermore, \cref{sn} is equivalent to
\eqe{
    x_i = \frac{1}{n}\sum_{j=1}^{n}(x_j-\nabla f_j(x_j)), \ i \in \mathcal{V}.
    \nonumber
}
Therefore, if each agent $i$ takes the following dynamics:
\eqe{ \label{gf}
    \dot{x}_i(t) &= \frac{1}{n}\sum_{j=1}^{n}(x_j-\nabla f_j(x_j)) - x_i(t) \\
    &= -\frac{1}{n}\sum_{j=1}^{n} \nabla f_j(x_j(t)) - \lt(x_i(t)-\frac{1}{n}\sum_{j=1}^{n}x_j(t)\rt),
}
% which implies that \cref{disp} will be solved if all agents can track to $\frac{1}{n}\sum_{j=1}^{n}(x_j-\nabla f_j(x_j))$. Given the above observation, 
% , which leads to the following distributed algorithm
obviously $x_i(t)$ will converge to $x^*$. Though \cref{gf} is not distributed, it is natural to apply the continuous-time DAC algorithm proposed in \cite{kia2015dynamic} to track $\frac{1}{n}\sum_{j=1}^{n}(x_j-\nabla f_j(x_j))$ distributedly:
\eqe{ \label{ap4}
    \dot{\mx} &= - \alpha(\mx-(\mx-\nabla f(\mx))) -\mz - \beta\mL\mx, \\
    \dot{\mz} &= \alpha\beta\mL\mx,
    \nonumber
}
which is exactly \cref{ap}. As a distributed version of \cref{gf}, $\alpha\nabla f_i(x_i) + z_i$ can be naturally seen as asgent $i$'s estimation of the global gradient $\frac{1}{n}\sum_{j=1}^{n} \nabla f_j(x_j(t))$. In other words, \cref{ap} can still be seen as a kind of gradient tracking algorithm, but with the implicit tracking approach, where $\alpha\nabla f(\mx) + \mz$ plays the role in tracking the global gradient. Compared with explicit gradient tracking algorithms, such as DIGing and ATC tracking, there is no explicit state variable to track the global gradient in \cref{ap}, which is the reason we call the tracking approach used in \cref{ap} the implicit tracking approach. 
 % which is exactly the new perspective proposed in \cite{li2022implicit},

As mentioned before, the implicit tracking approach is used in the algorithm design for distributed constraint-coupled optimization problems in \cite{li2022implicit}. Different from \cite{li2022implicit}, we focus on the classical set-up \cref{cenp} and develop a unified algorithmic framework for it, with the help of a brand-new understanding of existing distributed optimization algorithms that the implicit tracking mechanism offers. 

% With the help of the idea that employing DAC algorithms to track $\frac{1}{n}\sum_{j=1}^{n}(x_j-\nabla f_j(x_j))$ distributedly, we can easily develop a DT counterpart of \cref{ap}:
% \eqe{ \label{best}
%     \mx^{k+2} = (\mI+\mW)\mx^{k+1} - \mx^k - \alpha(\nabla f(\mx^{k+1})-\nabla f(\mx^{k})),
% }

Besides the continuous-time algorithm \cref{ap}, the implicit tracking mechanism exists generally in many discrete-time algorithms, like EXTRA, NIDS, and Exact diffusion. Though these algorithms do not have an explicit state variable to track the global gradient, it is feasible to transform them to explicit tracking forms. 
% which we call implicit tracking algorithms, we can still transform them to explicit tracking forms by unfolding their implicit tracking structure. 
Let us take EXTRA and Exact diffusion as CTA and ATC examples respectively. EXTRA updates as
\eqe{
	\mx^{k+2} &= 2\tW\mx\+ - \tW\mx^k - \alpha(\nabla f(\mx\+)-\nabla f(\mx^k)),
    \nonumber
}
where $\tW = \frac{\mI+\mW}{2}$. Define
\eqe{
	\mg\+ = -\tW(\mx\+-\mx^k) + \alpha(\nabla f(\mx\+)-\nabla f(\mx^k)),
    \nonumber
}
we have
\eqe{
	\mx\+ = \tW\mx^k - \mg^k
    \nonumber
}
and
\eqe{
	\mg\+ =& \tW\mx^k -\tW\lt(\tW\mx^k - \mg^k\rt) + \alpha(\nabla f(\mx\+)-\nabla f(\mx^k)) \\
	=&  \tW\mg^k + \lt(\mI-\tW\rt)\tW\mx^k + \alpha(\nabla f(\mx\+)-\nabla f(\mx^k)),
    \nonumber
}
then we can obtain the explicit tracking form of EXTRA:
\eqe{
	\mx\+ &= \tW\mx^k - \mg^k, \\
	\mg\+ &= \tW\mg^k + (\mI-\tW)\tW\mx^k+ \alpha(\nabla f(\mx\+)-\nabla f(\mx^k)).
    \nonumber
}
Exact diffusion updates as
\eqe{
	\mx^{k+2} = \tW(2\mx\+-\mx^k-\alpha(\nabla f(\mx\+)-\nabla f(\mx^k))).
    \nonumber
}
Define
\eqe{
	\mg\+ = \mx^k-\mx\+ +\alpha(\nabla f(\mx\+)-\nabla f(\mx^k)),
    \nonumber
}
we have
\eqe{
	\mx\+ = \tW(\mx^k-\mg^k)
    \nonumber
}
and
\eqe{
	\mg\+ &= \mx^k-\tW(\mx^k-\mg^k) +\alpha(\nabla f(\mx\+)-\nabla f(\mx^k)) \\
	&= \tW\mg^k + (\mI-\tW)\mx^k+\alpha(\nabla f(\mx\+)-\nabla f(\mx^k)),
    \nonumber
}
then the explicit tracking form of Exact diffusion is given as
\eqe{
	\mx\+ &= \tW(\mx^k-\mg^k), \\
	\mg\+ &=\tW\mg^k + (\mI-\tW)\mx^k+\alpha(\nabla f(\mx\+)-\nabla f(\mx^k)).
    \nonumber
}

Apart from EXTRA and Exact diffusion, many other distributed optimization algorithms also have the implicit gradient tracking structure and we can transform them to explicit tracking forms, which inspires the design of the unified framework. Our unified framework, i.e., UGT, has two variants, one is CTA-UGT:
\eqe{\label{CTA-UGT}
    \mx\+ &= \mW_1\mx^k - \mg^k, \\
    \mg\+ &= \mW_2\mg^k + \beta(\mI-\mW_2)\mW_1\mx^k \\
    &\quad + \alpha(\nabla f(\mx\+)-\nabla f(\mx^k)),
    \nonumber
}
and another one is ATC-UGT:
\eqe{\label{ATC-UGT}
    \mx\+ &= \mW_1(\mx^k - \mg^k), \\
    \mg\+ &= \mW_2\mg^k + \beta(\mI-\mW_2)\mx^k \\
    &\quad + \alpha(\nabla f(\mx\+)-\nabla f(\mx^k)),
    \nonumber
}
where $\beta\geq0$ is a tunable parameter, $\mW_1 = W_1 \otimes I_m$, $\mW_2 = W_2 \otimes I_m$, $W_1$ and $W_2$ are weight matrices, and $\mg^0 = \alpha \nabla f(\mx^0)$. By choosing different $\mW_1$, $\mW_2$, and $\beta$, UGT can recover most existing algorithms, as shown in \cref{tab1}.
\begin{remark}
As a matter of fact, UGT is a generalized gradient tracking algorithm. Compared with the two classical gradient tracking algorithms (DIGing and ATC tracking), UGT has an extra term: for CTA-UGT, it is $\beta(\mI-\mW_2)\mW_1\mx^k$; for ATC-UGT, it is $\beta(\mI-\mW_2)\mx^k$, which we call the modified term. It is worth noting that the introduce of the modified term does not increase the number of communication rounds in each iteration, which is obvious for ATC-UGT, while for CTA-UGT, we only need to notice that the updating of $\mg^k$ can be rewritten as
\eqe{
    \mg\+ &= (1-\beta)\mW_2\mg^k + \beta\lt[\mg^k+(\mI-\mW_2)\mx\+\rt]  \\
    &\quad + \alpha(\nabla f(\mx\+)-\nabla f(\mx^k)),
    \nonumber
}
hence the number of communication rounds of CTA-UGT in each iteration is still $2$.
\end{remark}

\begin{remark}
The parameter $\beta$ in the modified term plays a critical role in the performance of UGT. In existing algorithms, there are only two possible values of $\beta$: $0$ and $1$, as shown in \cref{tab1}. Nevertheless, $\beta$ can be set to any values, as long as the convergence can be guaranteed (the upper bound of $\beta$ is derived in \cref{convergence}). In our opinion, UGT with different $\beta$ are essentially different algorithms, and whose convergence rates may have huge differences in numerical experiments. Therefore, it is convenient to develop efficient versions of UGT by choosing an appropriate value for $\beta$. 
% In numerical experiments, it is found that we can obtain a better convergence rate if $\beta$ is appropriately selected.
\end{remark}
% To explain the role it plays in UGT, let us first consider CTA-UGT. In CTA-UGT, $\mg\+$ can be rewritten as
% \eqe{
% 	\mg\+ &= (1-\beta)\mW_2\mg^k + \beta\lt[\mg^k+(\mI-\mW_2)\mx\+\rt]  \\
%     &\quad + \alpha(\nabla f(\mx\+)-\nabla f(\mx^k)),
% }
% where $(1-\beta)\mW_2\mg^k$ and $\beta\lt[\mg^k+(\mI-\mW_2)\mx\+\rt]$ represent two different strategies that make $\mg^k$ achieve consensus, as a matter of fact. The first is the classical one, which is easy to understand, but the second is not very straightforward. From \cref{DIGing}, we can see that $\mx^k$ and $\mg^k$ influences each other,  and one can not converge unless the other one converges, meanwhile, one will converge if the other one converges. Considering of that, the idea behind $\beta\lt[\mg^k+(\mI-\mW_2)\mx\+\rt]$ is natural: besides the direct strategy, to make $\mg^k$ achieve consensus, we can also adopt the indirect strategy, i.e, making $\mg^k$ change along the direction which is in favor of making $\mx^k$ achieve consensus, in turn, that will drive $\mg^k$ to achieve consensus. As shown in \cref{tab1}, existing algorithms either choose the direct strategy (i.e., $\beta=0$), or choose the indirect strategy (i.e., $\beta=1$). It is natural to think that if we can obtain faster algorithms when we mix the direct strategy and the indirect strategy together (i.e., $0<\beta<1$), which is verified in numerical experiments.
% If we further assume the symmetry of $W_1$ and $W_2$, then we can also obtain a primal-dual explaination of $\beta\lt[\mg^k+(\mI-\mW_2)\mx\+\rt]$, but 

% \begin{remark}

% \end{remark}

\begin{figure*}[tb]
\begin{center}
    \subfigure[Directed circle graph]{
    \begin{minipage}[b]{0.4\textwidth}
        \includegraphics[scale=0.35]{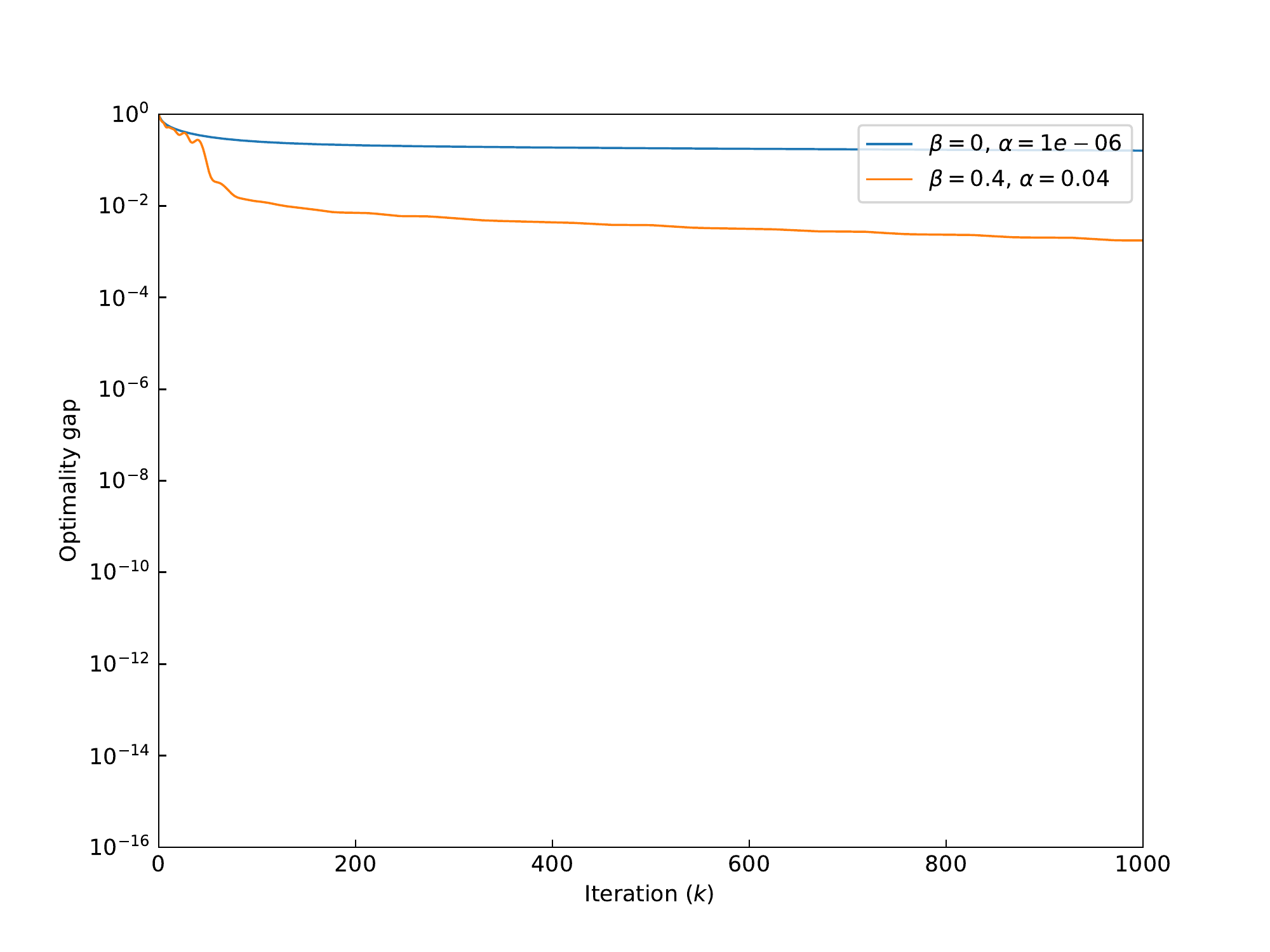}
    \end{minipage}
    }
    \subfigure[Directed exponential graph with $e=2$]{
    \begin{minipage}[b]{0.4\textwidth}
        \includegraphics[scale=0.35]{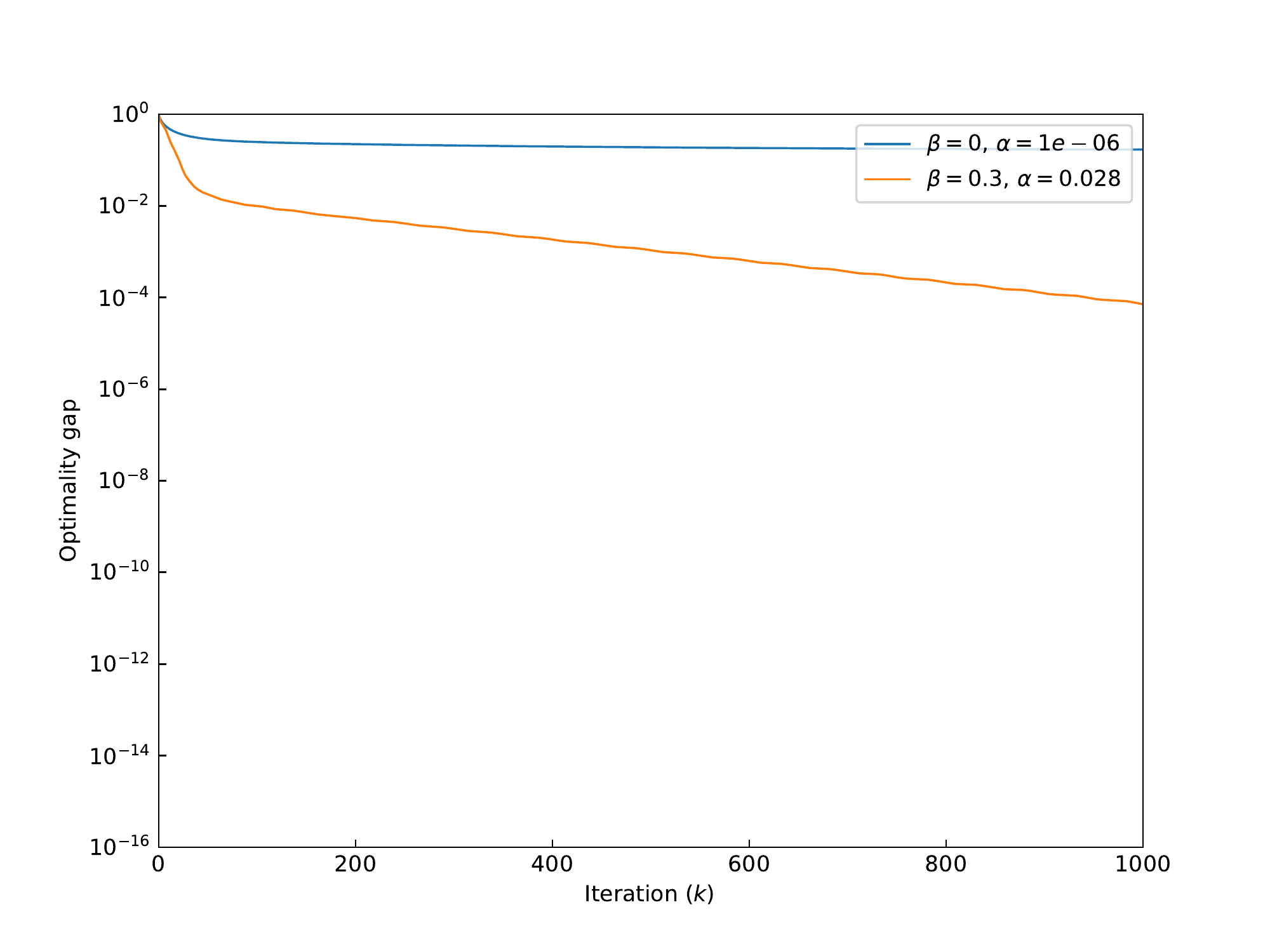}
    \end{minipage}
    }
    \subfigure[Directed exponential graph with $e=4$]{
    \begin{minipage}[b]{0.4\textwidth}
        \includegraphics[scale=0.35]{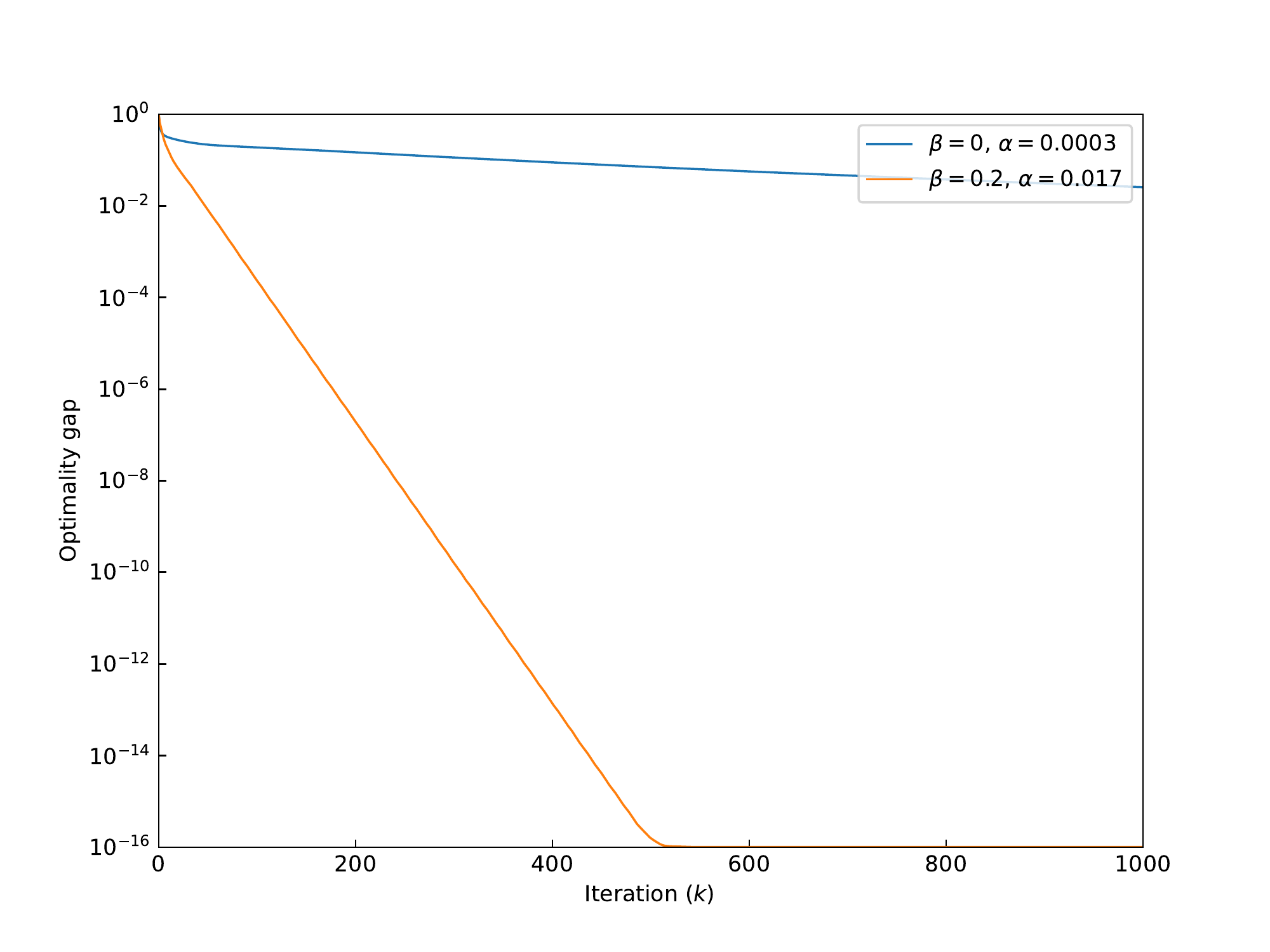}
    \end{minipage}
    }
    \subfigure[Directed exponential graph with $e=6$]{
    \begin{minipage}[b]{0.4\textwidth}
        \includegraphics[scale=0.35]{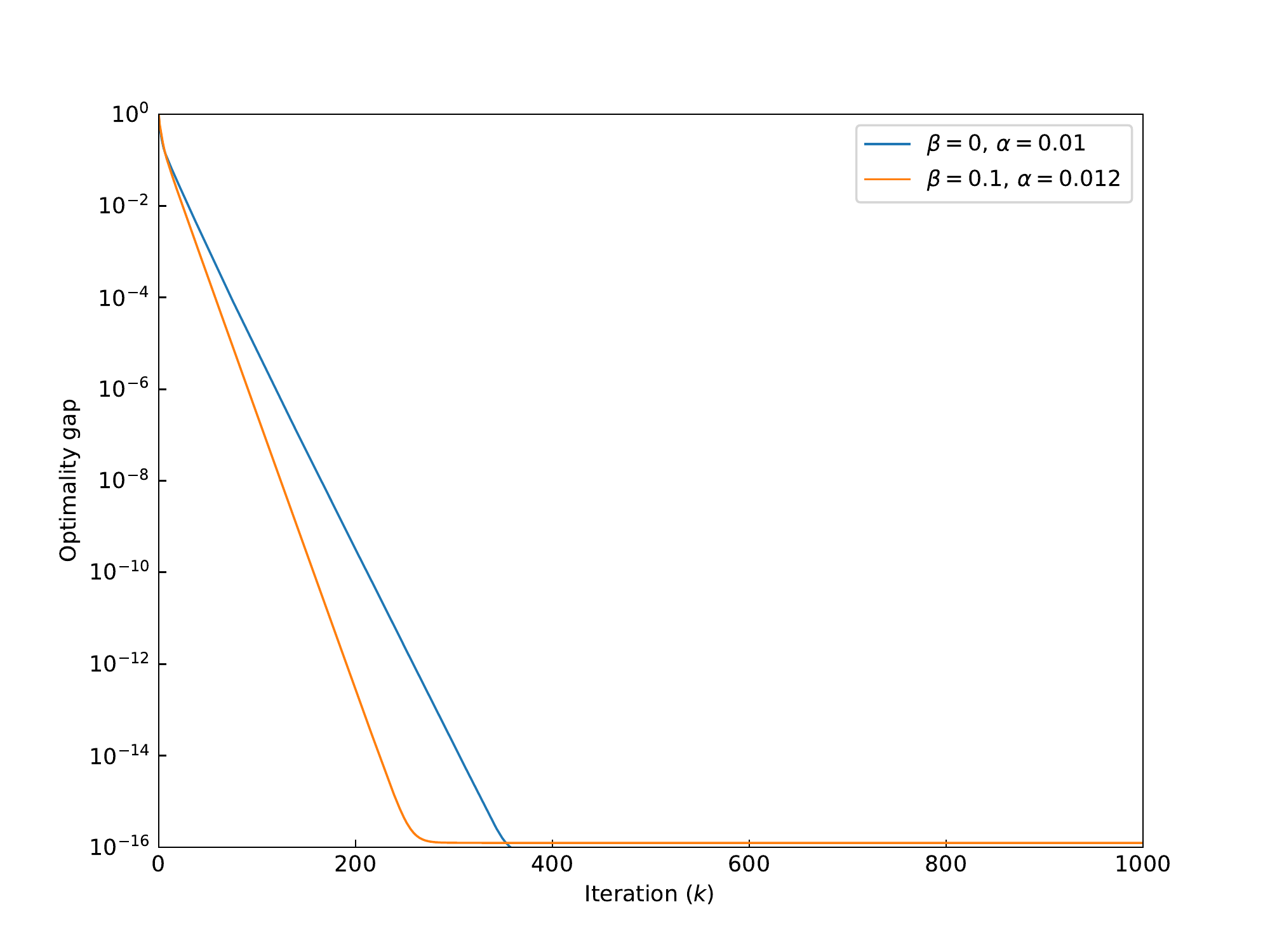}
    \end{minipage}
    }
\caption {The convergence rates of different versions of CTA-UGT over digraphs.}
\label{fig1}
\end{center}
\end{figure*}

\begin{figure*}[tb]
\begin{center}
    \subfigure[Directed circle graph]{
    \begin{minipage}[b]{0.4\textwidth}
        \includegraphics[scale=0.35]{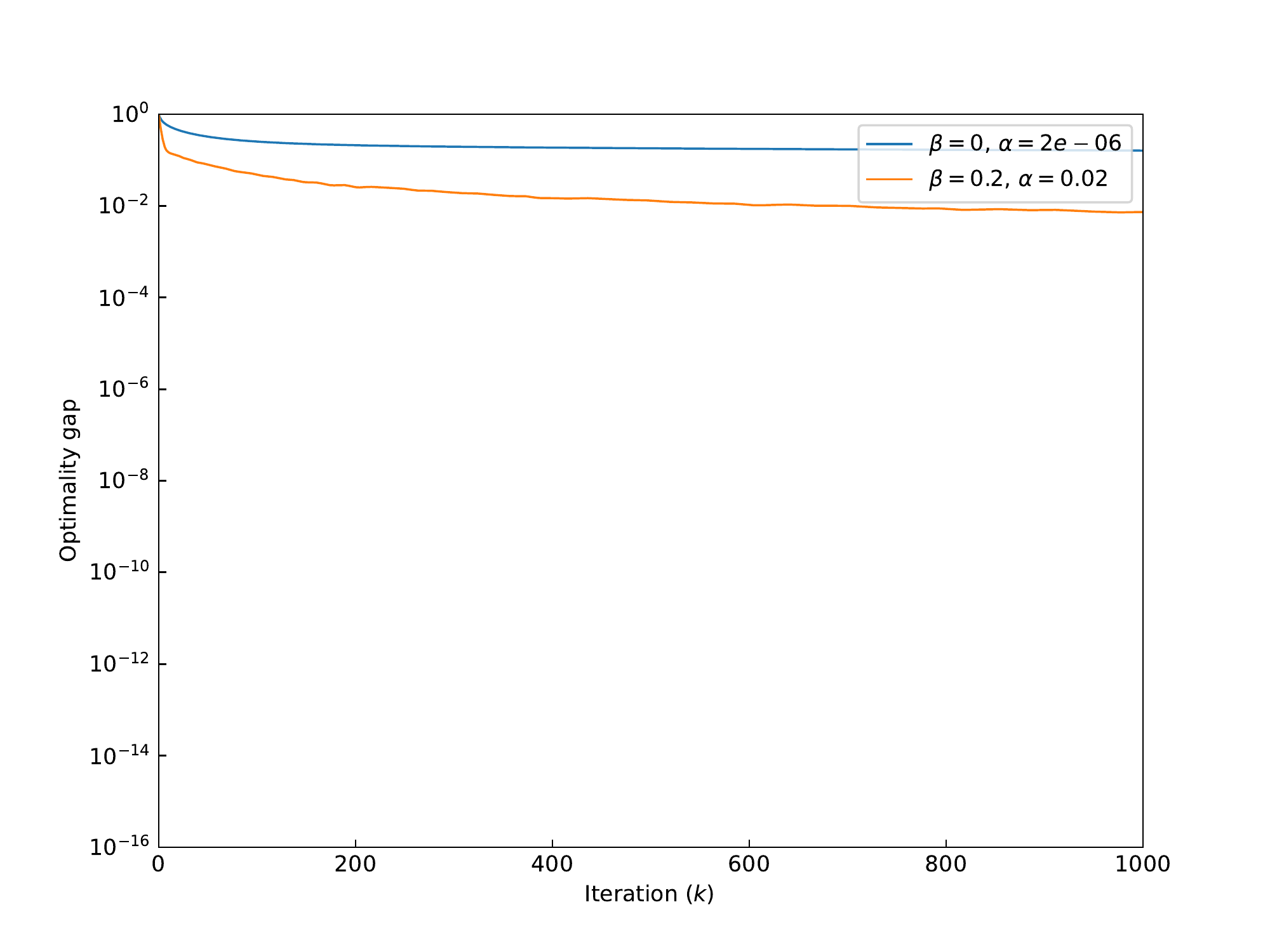}
    \end{minipage}
    }
    \subfigure[Directed exponential graph with $e=2$]{
    \begin{minipage}[b]{0.4\textwidth}
        \includegraphics[scale=0.35]{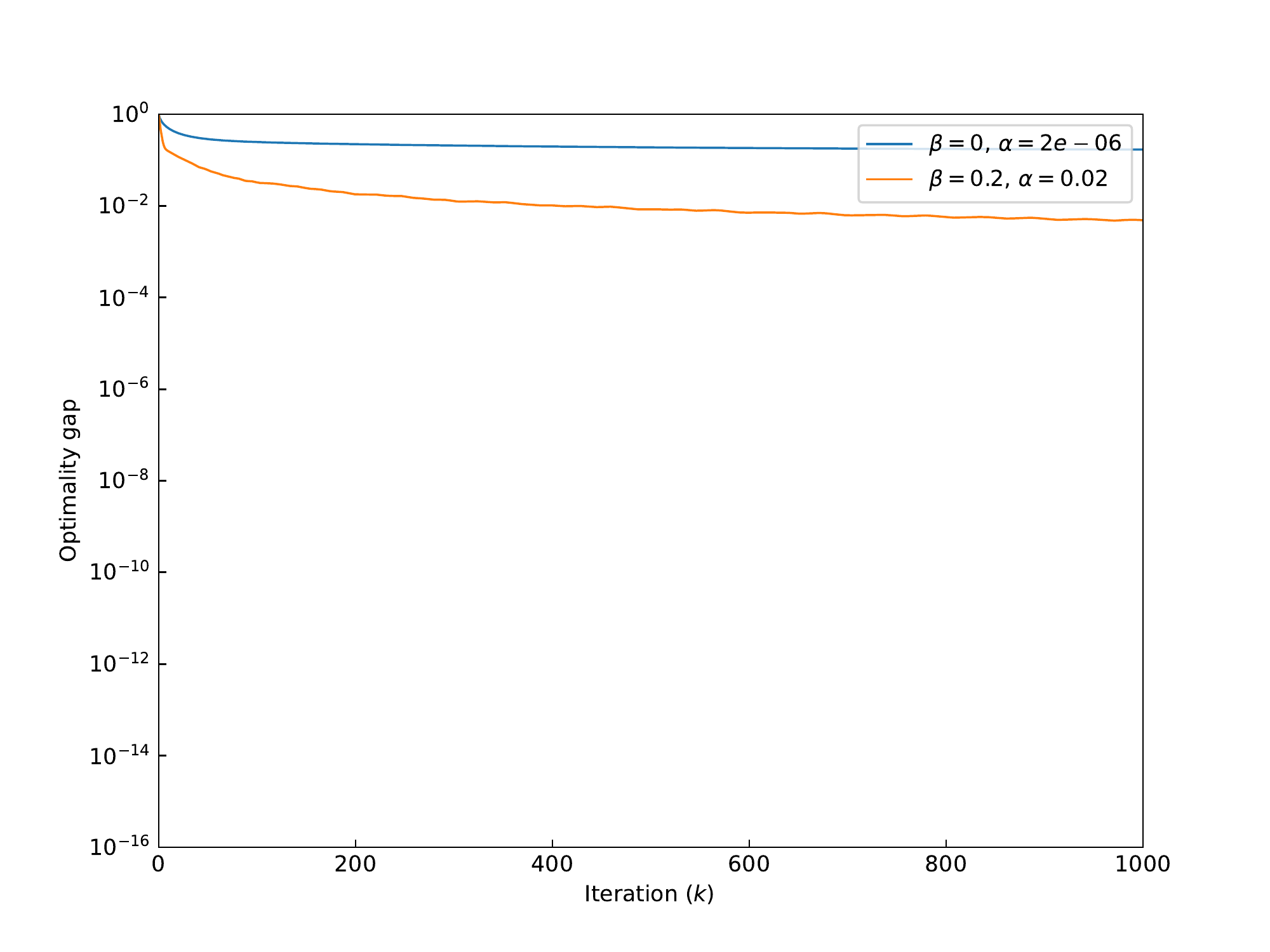}
    \end{minipage}
    }
    \subfigure[Directed exponential graph with $e=4$]{
    \begin{minipage}[b]{0.4\textwidth}
        \includegraphics[scale=0.35]{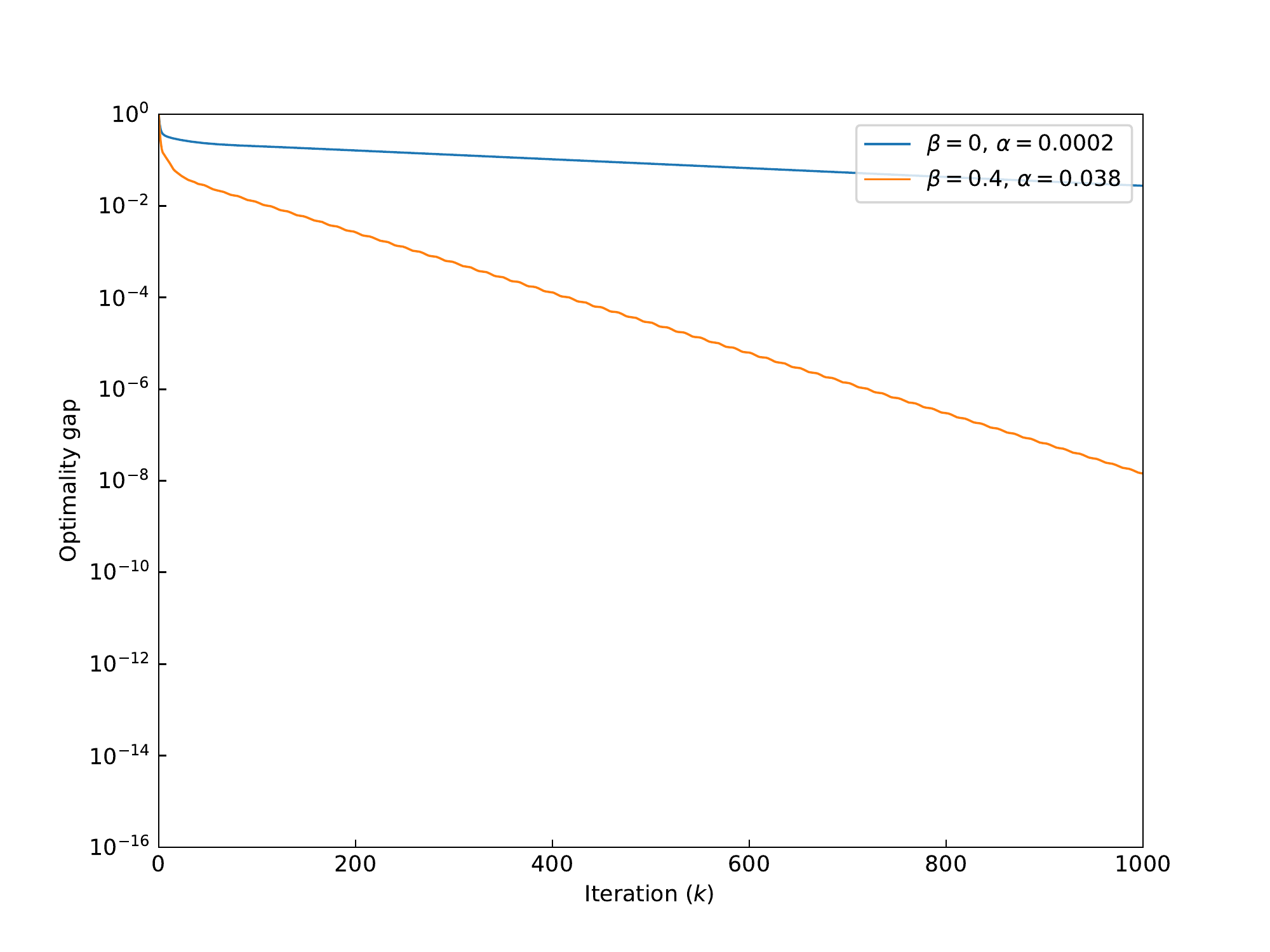}
    \end{minipage}
    }
    \subfigure[Directed exponential graph with $e=6$]{
    \begin{minipage}[b]{0.4\textwidth}
        \includegraphics[scale=0.35]{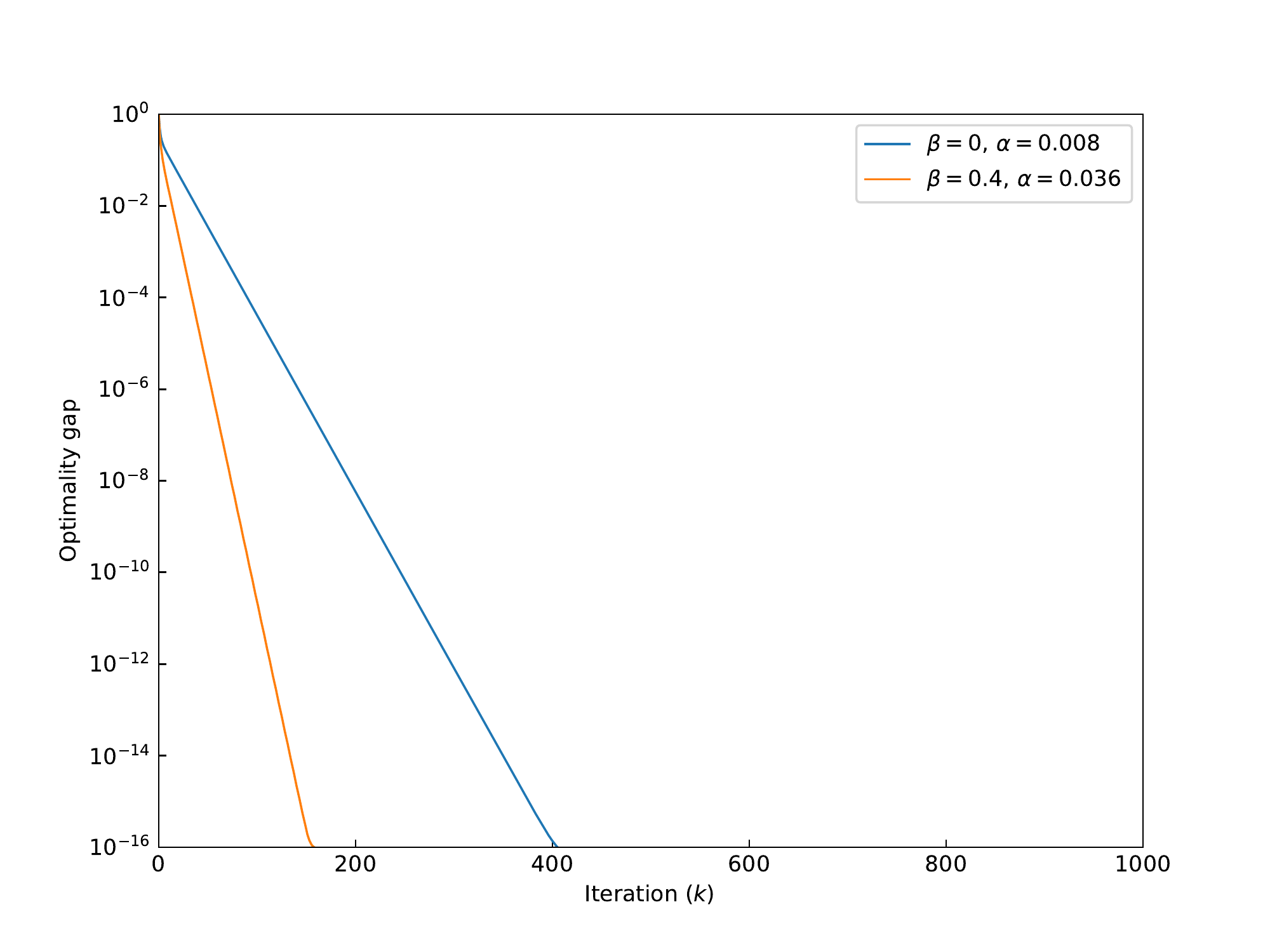}
    \end{minipage}
    }
\caption {The convergence rates of different versions of ATC-UGT over digraphs.}
\label{fig2}
\end{center}
\end{figure*}

\begin{figure*}[tb]
\begin{center}
    \subfigure[Circle graph]{
    \begin{minipage}[b]{0.4\textwidth}
        \includegraphics[scale=0.35]{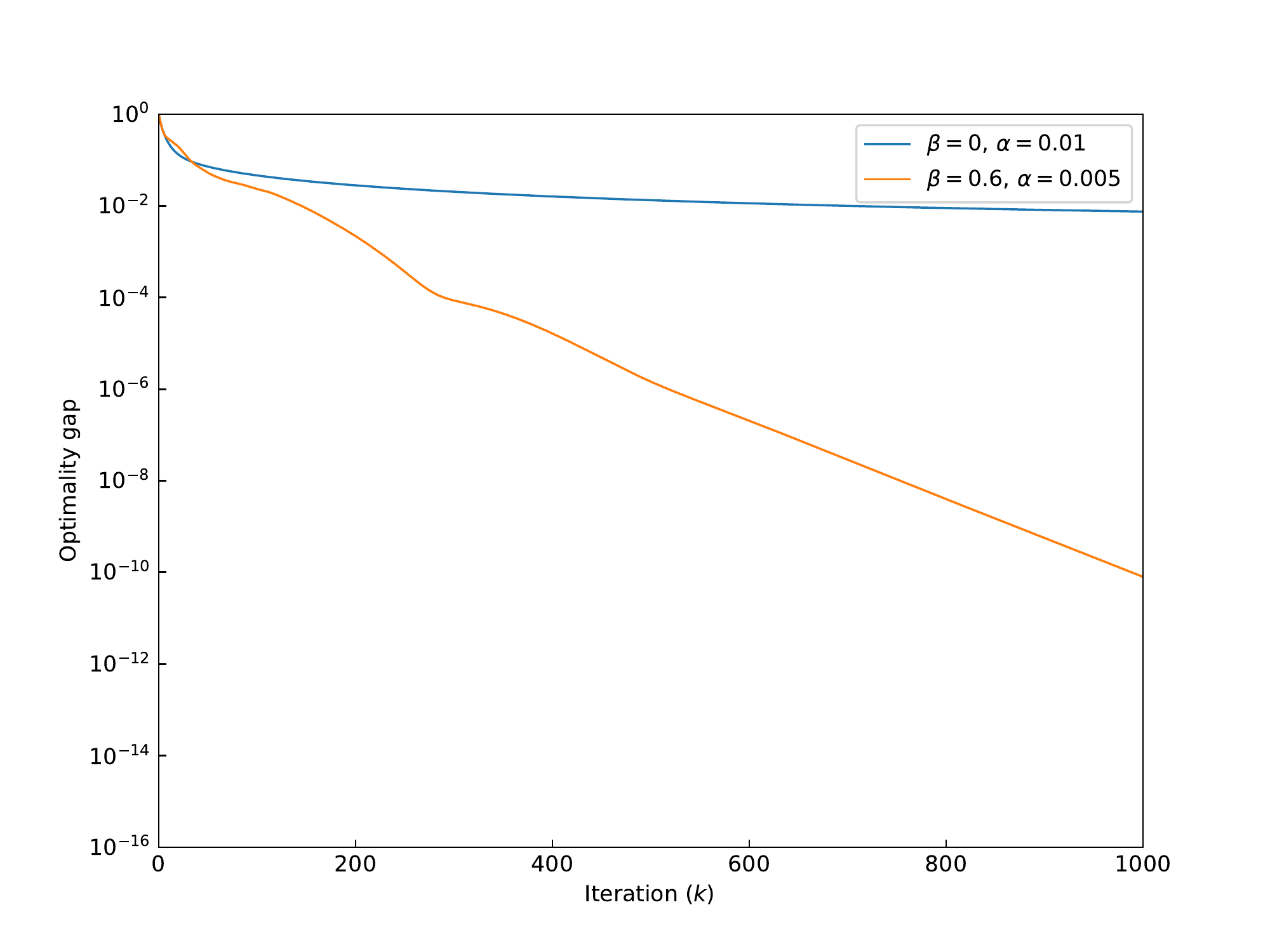}
    \end{minipage}
    }
    \subfigure[Random graph with $p=0.05$]{
    \begin{minipage}[b]{0.4\textwidth}
        \includegraphics[scale=0.35]{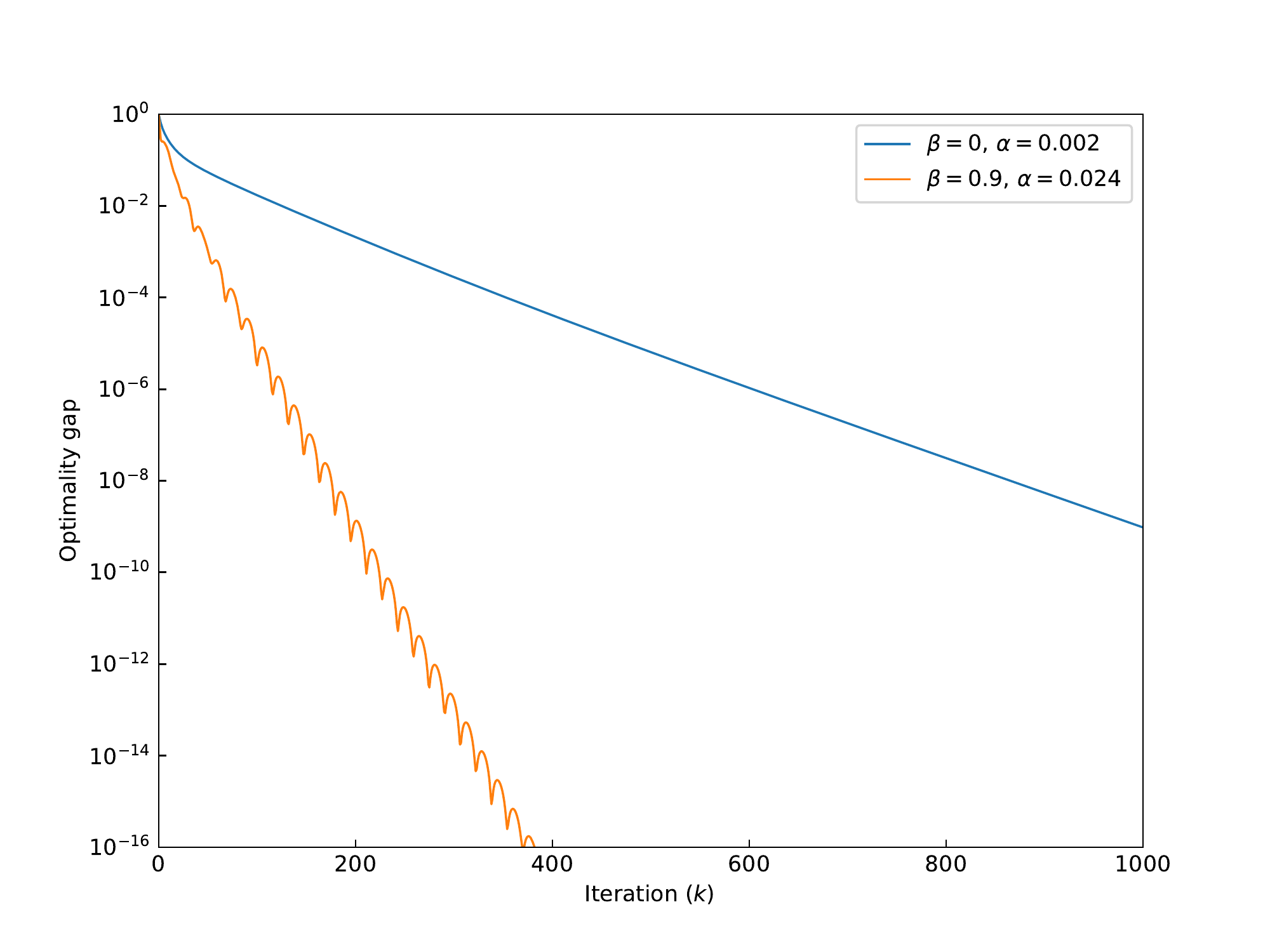}
    \end{minipage}
    }
    \subfigure[Random graph with $p=0.1$]{
    \begin{minipage}[b]{0.4\textwidth}
        \includegraphics[scale=0.35]{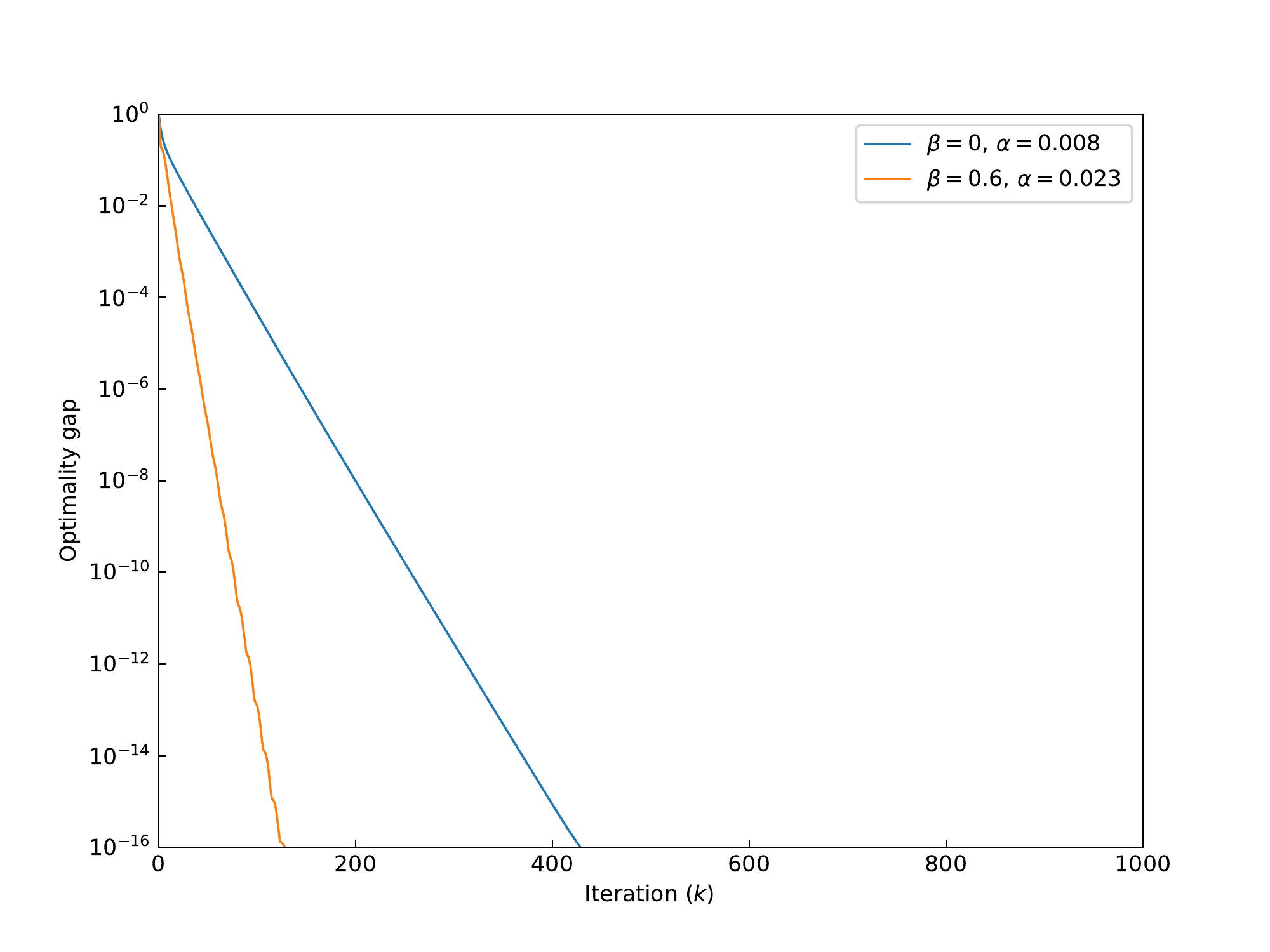}
    \end{minipage}
    }
    \subfigure[Random graph with $p=0.3$]{
    \begin{minipage}[b]{0.4\textwidth}
        \includegraphics[scale=0.35]{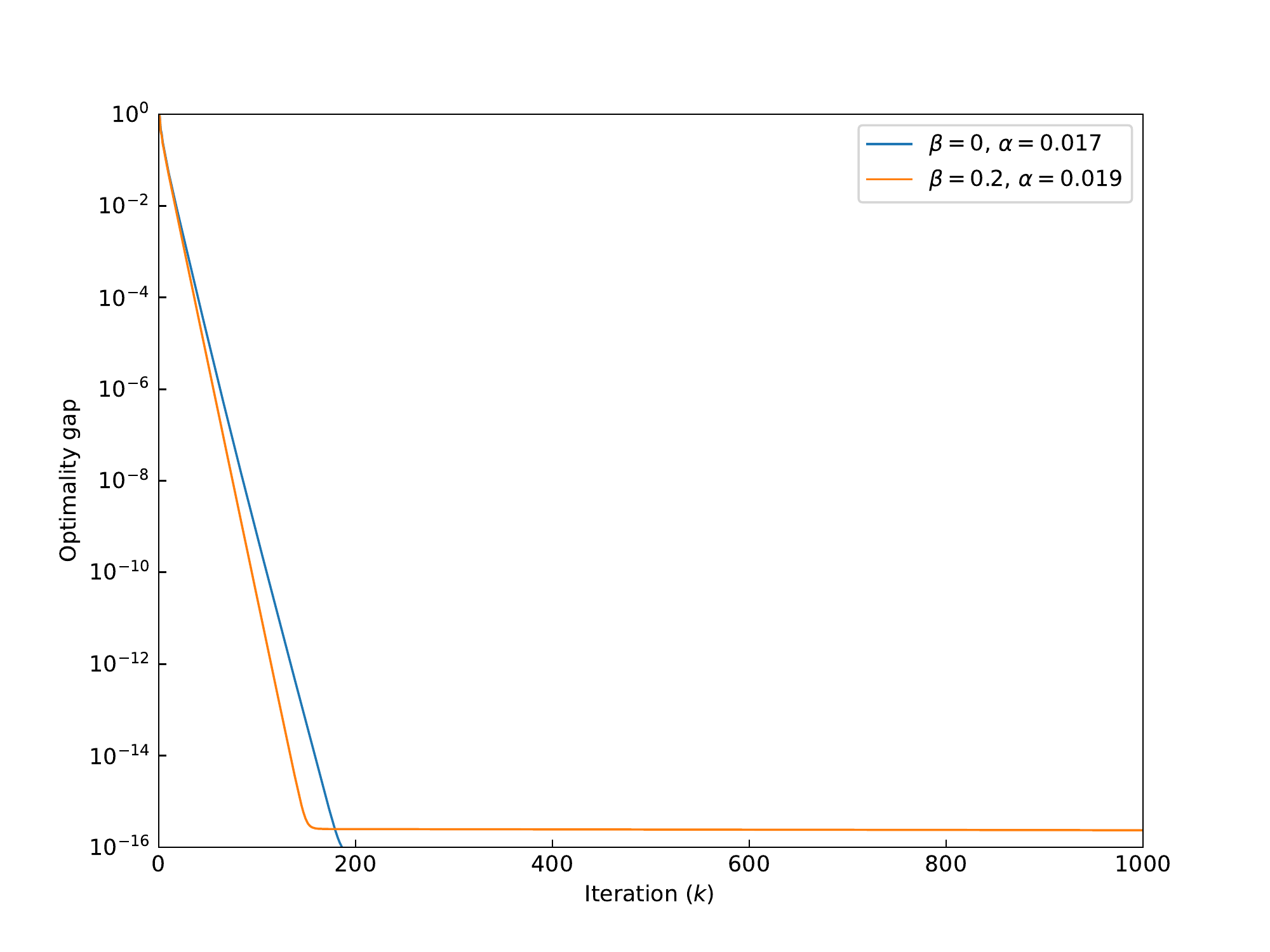}
    \end{minipage}
    }
\caption {The convergence rates of different versions of CTA-UGT over undigraphs.}
\label{fig3}
\end{center}
\end{figure*}

\begin{figure*}[tb]
\begin{center}
    \subfigure[Circle graph]{
    \begin{minipage}[b]{0.4\textwidth}
        \includegraphics[scale=0.35]{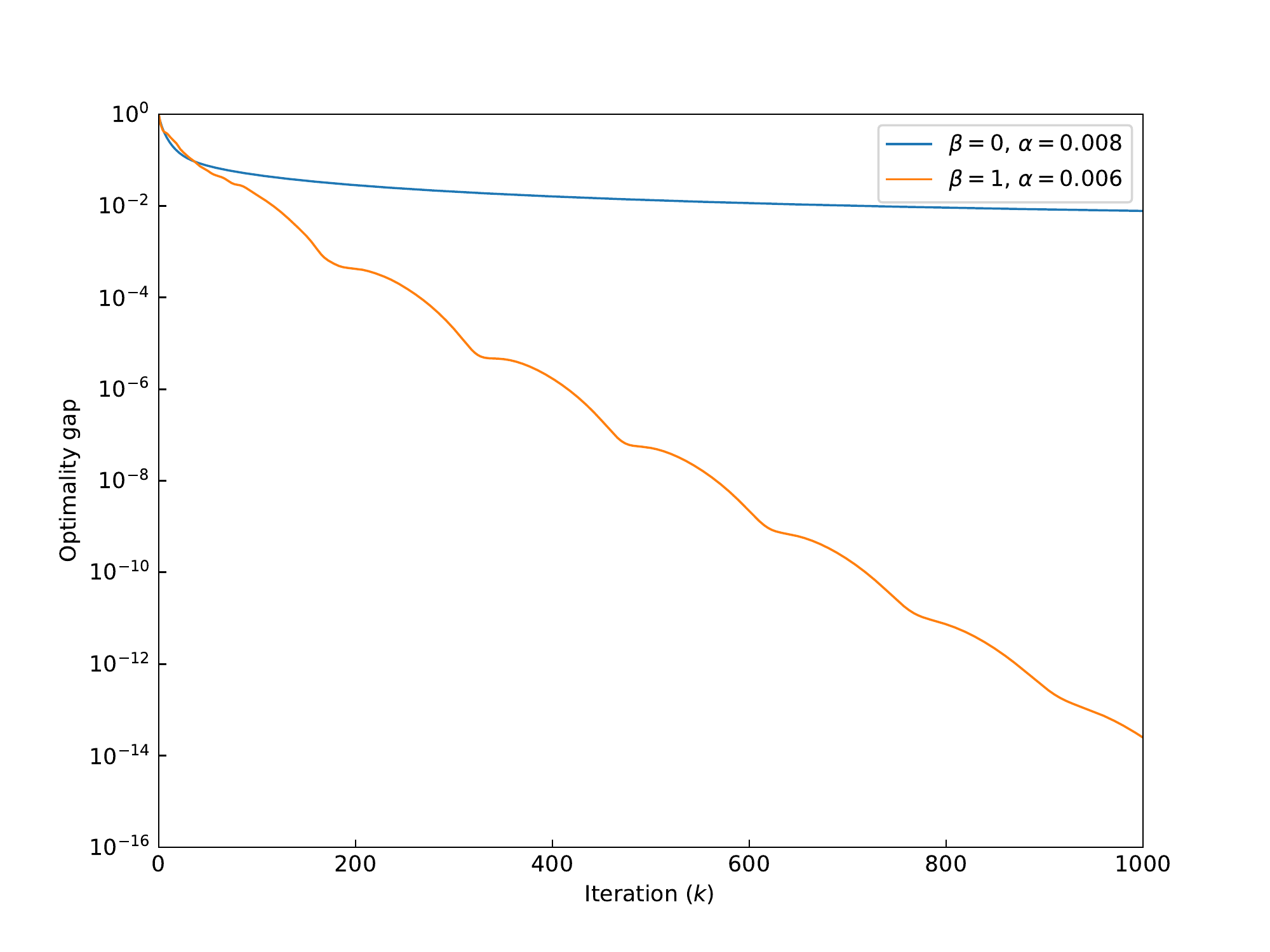}
    \end{minipage}
    }
    \subfigure[Random graph with $p=0.05$]{
    \begin{minipage}[b]{0.4\textwidth}
        \includegraphics[scale=0.35]{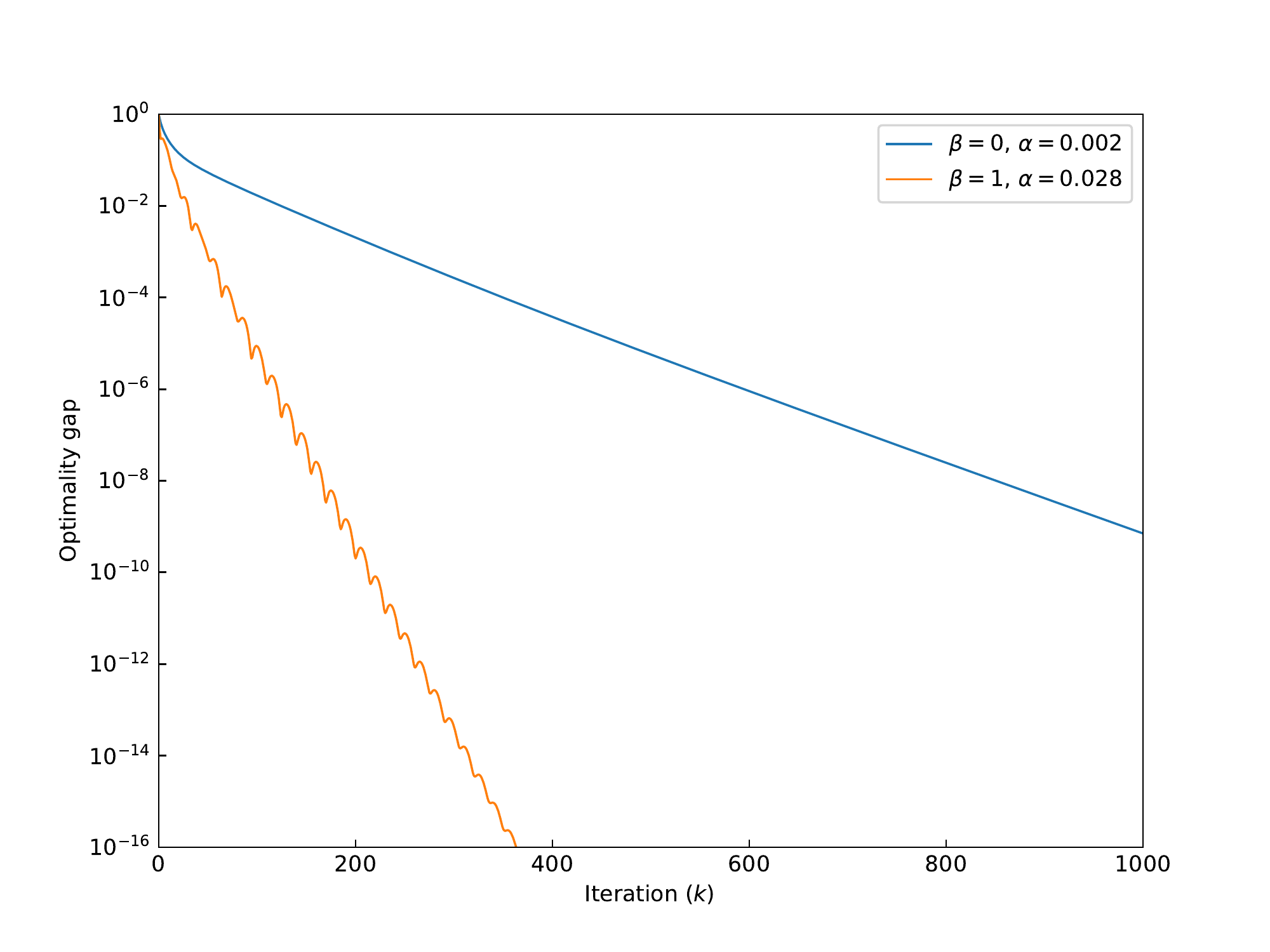}
    \end{minipage}
    }
    \subfigure[Random graph with $p=0.1$]{
    \begin{minipage}[b]{0.4\textwidth}
        \includegraphics[scale=0.35]{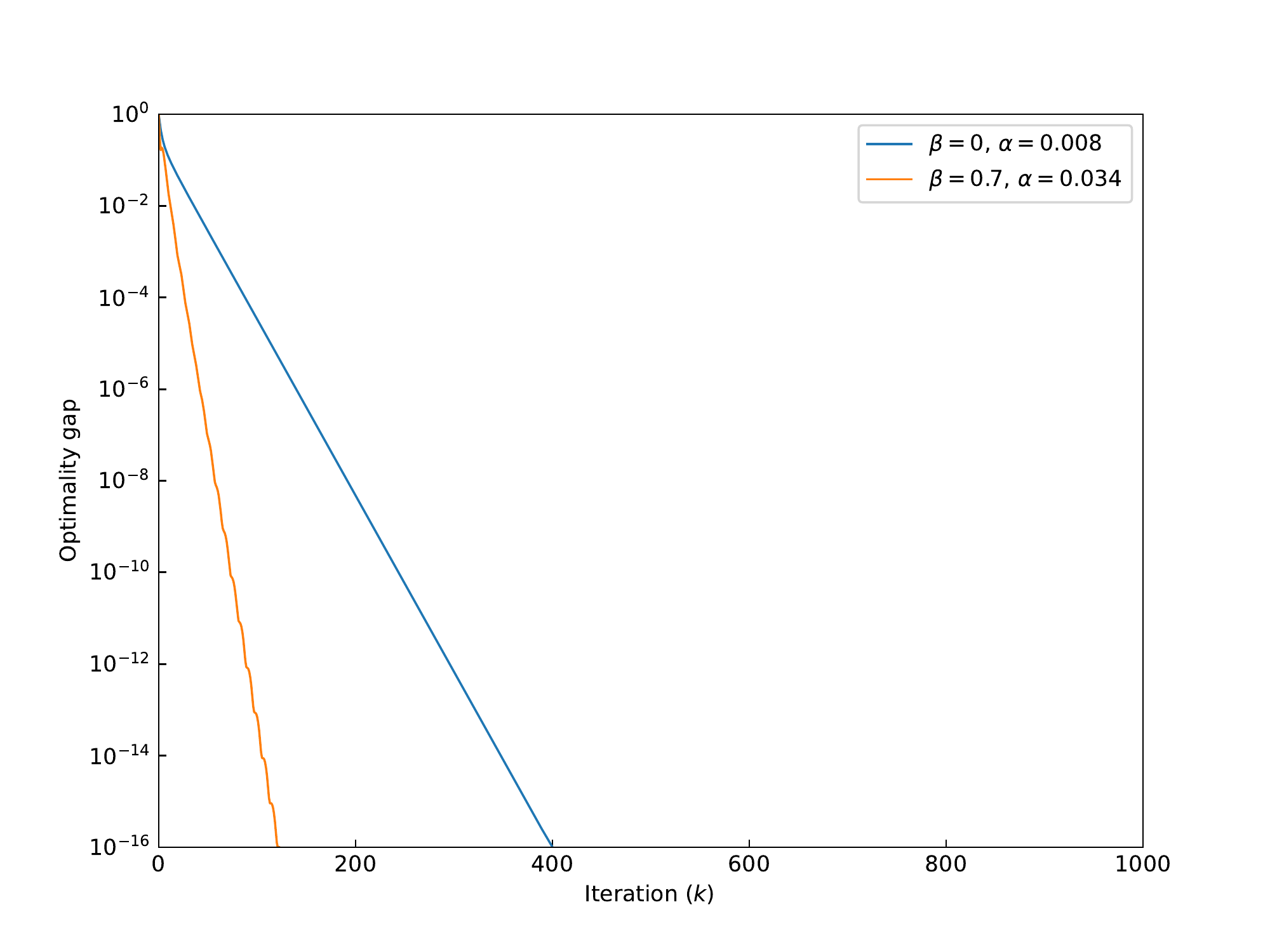}
    \end{minipage}
    }
    \subfigure[Random graph with $p=0.3$]{
    \begin{minipage}[b]{0.4\textwidth}
        \includegraphics[scale=0.35]{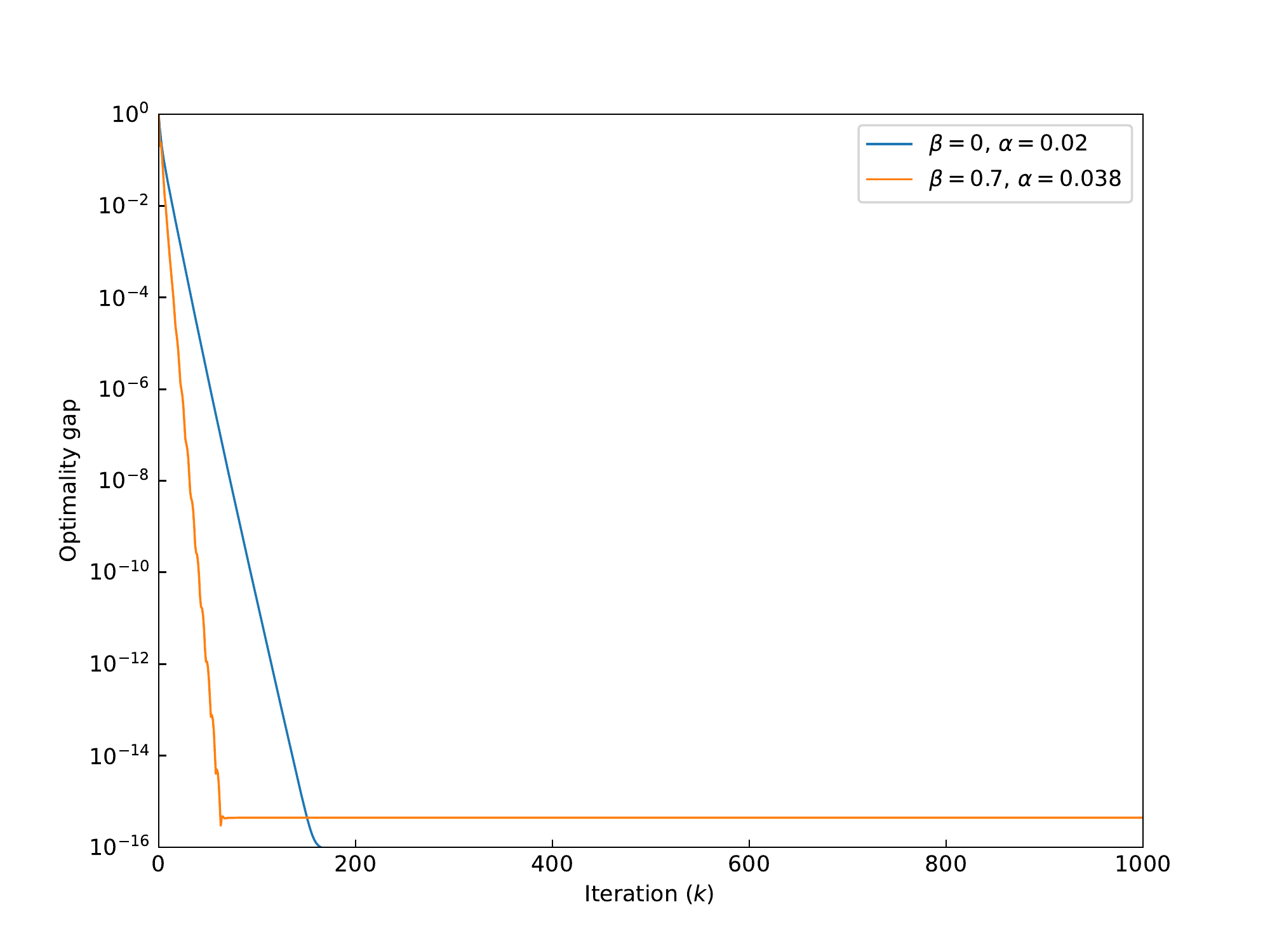}
    \end{minipage}
    }
\caption {The convergence rates of different versions of ATC-UGT over undigraphs.}
\label{fig4}
\end{center}
\end{figure*}

\section{Convergence Analysis} \label{convergence}
% In this section, we prove the linear convergence of CTA-UGT and ATC-UGT under the following assumptions.
We first give the following assumption about weight matrices $W_1$ and $W_2$. 
\begin{assumption} \label{W}
$W_1$ and $W_2$ are primitive and doubly stochastic.
\end{assumption}

$W_1$ and $W_2$ are associated with $\mathcal{G}$. $W_1(W_2) = [w_{ij}]$ can be constructed as: $w_{ii}>0$; $w_{ij}>0$ if $(i, j) \in \mathcal{E}$, otherwise $w_{ij}=0$. Based on this constructed rule, $W_1$ and $W_2$ are primitive if $\mathcal{G}$ is strongly connected.
Given \cref{W}, $W_1$ and $W_2$ have the following property \cite{qu2017harnessing}:
\eqe{
    \sigma_1 =& \left\|W_1-\frac{1}{n}\1\1^{\top}\right\| \in (0,1), \\
    \sigma_2 =& \left\|W_2-\frac{1}{n}\1\1^{\top}\right\| \in (0,1).
    \nonumber
}
Though the above property is derived under the assumption that $\mathcal{G}$ is undirected and connected in \cite{qu2017harnessing}, it is trivial to prove it for our case by feat of the Perron–Frobenius theory.
% , hence we omit the proof.
\begin{remark}
Different from \cite{alghunaim2020decentralized} and \cite{xu2021distributed}, we do not assume the symmetry of $W$. Consequently, \cref{W} can not only be satisfied by undigraphs, but also digraphs that permit doubly stochastic weight matrices, which is the reason we use the digraph to model the communication topology. On the converse, the communication topology can only be modeled by an undigraph in \cite{alghunaim2020decentralized} and \cite{xu2021distributed}, due to the symmetry of $W$. 
\end{remark}

\begin{assumption} \label{str}
$F$ is $\mu$-strongly convex and $f_i$ is $l_i$-smooth, $i \in \mathcal{V}$.
\end{assumption}

Let $l = \min_{i \in \mathcal{V}}l_i$, $\bar{l}=\frac{1}{n}\sum_{i=1}^nl_i$, $\eta = \min\{|1-\alpha \mu|, |1-\alpha \bar{l}|\}$, $\sigma_1 \in (0,1)$ and $\sigma_2 \in (0,1)$ be the second largest singular values of $W_1$ and $W_2$ respectively. Define two non-negative matrices as
$$G_1(\alpha, \beta) = \left[\begin{array}{cccc}
\sigma_1 & 1 & 0 \\
\Centerstack{$\beta\sigma_1(1+\sigma_2)$ \\ $+ \alpha l(1+\sigma_1+\alpha l)$} & \sigma_2+\alpha l & (\alpha l)^2 \\
\alpha l & 0 & \eta
\end{array}\right] $$
and
$$G_2(\alpha, \beta) = \left[\begin{array}{cccc}
\sigma_1 & \sigma_1 & 0 \\
\Centerstack{$\beta(1+\sigma_2)$ \\ $+ \alpha l(1+\sigma_1+\alpha l)$} & \sigma_2+\alpha l & (\alpha l)^2 \\
\alpha l & 0 & \eta
\end{array}\right]. $$
\begin{theorem} \label{the1}
Suppose \cref{W,str} hold and choose $\alpha>0$ and $\beta\geq0$ such that $\rho(G_1(\alpha, \beta)) < 1$, then $\mx^k$ generated by CTA-UGT converges to $\mx^*$ with the R-linear rate $O(\rho(G_1(\alpha, \beta)^k)$.
Furthermore, $\rho(G_1(\alpha, \beta)) < 1$ if
$$0\leq\beta<\frac{(1-\sigma_1)(1-\sigma_2)}{\sigma_1(1+\sigma_2)}$$
and
$$0<\alpha<\min\left\{\frac{1-\sigma_2}{l}, \frac{-b+\sqrt{b^2-4ac}}{2a}\right\},$$
where $a=l^3+\mu l^2$, $b=2\mu l$, and $c=\mu[\beta\sigma_1(1+\sigma_2)-(1-\sigma_1)(1-\sigma_2)]$.
\end{theorem}

\begin{theorem} \label{the2}
Suppose \cref{W,str} hold and choose $\alpha>0$ and $\beta\geq0$ such that $\rho(G_2(\alpha, \beta)) < 1$, then $\mx^k$ generated by ATC-UGT converges to $\mx^*$ with the R-linear rate $O(\rho(G_2(\alpha, \beta)^k)$.
Furthermore, $\rho(G_2(\alpha, \beta)) < 1$ if
$$0\leq\beta<\frac{(1-\sigma_1)(1-\sigma_2)}{\sigma_1(1+\sigma_2)}$$
and
$$0<\alpha<\min\left\{\frac{1-\sigma_2}{l}, \frac{-b+\sqrt{b^2-4ac}}{2a}\right\},$$
where $a=\sigma_1(l^3+\mu l^2)$, $b=(1+\sigma_1^2)\mu l$, and $c=\mu[\beta\sigma_1(1+\sigma_2)-(1-\sigma_1)(1-\sigma_2)]$.
\end{theorem}

For the sake of readability, the proofs of \cref{the1,the2} are placed in the appendix.
% \begin{remark}
% According to \cref{the1,the2}, the upper bound of $\beta$ depends on the 
% \end{remark}

\section{Numerical Experiments} \label{simu}
\newcolumntype{Y}{>{\centering\arraybackslash}X}
\strutlongstacks{T}
In this section, we evaluate the performance of UGT by solving the following quadratic programming problem:
\eqe{
     \min_{x \in \mathbb{R}^m} \ &F(x) = \frac{1}{n}\sum_{i=1}^n x_i^{\top}Q_ix_i + c_i^{\top}x_i, \\
     \nonumber
}
where $Q_i \in \mathbb{R}^{m \times m}$ is a positive definite matrix. In particular, different kinds of graphs are considered, included four digraphs: directed cycle graph, directed exponential graphs with $e = 2, 4, 6$, and four undigraphs: cycle graph, Erdos–Renyi random graphs with the connectivity probability $p = 0.05, 0.1, 0.3$. The directed exponential graph with $n$ nodes is generated by the rule that each node $i$ can send information to the nodes $(i+2^j)\mod n$, $j=0, 1, \cdots, e-1$. The weight matrices $W_1$ and $W_2$ are constructed by the Laplacian method \cite{shi2015extra}.

In experiments, we set $m = 2$ and $n = 100$, then generate $Q_i$ and $c_i$ randomly. For each graph, we compare the convergence rates of two versions of CTA-UGT (ATC-UGT): $\beta = 0$ (which corresponds to traditional gradient tracking algorithms) and $\beta = w$, where $w \geq 0$ is the constant chosen according to the rule that maximizing the convergence rate of CTA-UGT (ATC-UGT) for the given graph. For a given $\beta$, $\alpha$ is determined by the same rule. The experiments results are presented in \cref{fig1,fig2,fig3,fig4}, where the relative optimality gap is defined as $\frac{\|\mx(k)-\mx^*\|}{\|\mx(0)-\mx^*\|}$. It is easy to see that the latter usually has a far faster convergence rate, especially when the graph is poorly connected.

% \begin{figure*}[tb]
% \begin{center}
%     \subfigure[Circle graph, $\eta_2(L) \approx 0.02$]{
%     \begin{minipage}[b]{0.23\textwidth}
%         \includegraphics[scale=0.22]{figs/CTA-d-c.pdf}
%     \end{minipage}
%     }
%     \subfigure[Random graph with $p=0.05$, $\eta_2(L) \approx 0.21$]{
%     \begin{minipage}[b]{0.23\textwidth}
%         \includegraphics[scale=0.22]{figs/CTA-d-e-2.pdf}
%     \end{minipage}
%     }
%     \subfigure[Random graph with $p=0.1$, $\eta_2(L) \approx 0.54$]{
%     \begin{minipage}[b]{0.23\textwidth}
%         \includegraphics[scale=0.22]{figs/CTA-d-e-4.pdf}
%     \end{minipage}
%     }
%     \subfigure[Random graph with $p=0.3$, $\eta_2(L) \approx 6.34$]{
%     \begin{minipage}[b]{0.23\textwidth}
%         \includegraphics[scale=0.22]{figs/CTA-d-e-6.pdf}
%     \end{minipage}
%     }
% \caption {Experiment results for Case 1.}
% \label{fig1}
% \end{center}
% \end{figure*}

\section{Conclusion}
\label{conclusion}
In this work, the classical distributed optimization problem is studied. Inspired by the implicit tracking mechanism, UGT, a unified algorithmic framework, is developed from a pure primal perspective, which can unify most existing distributed optimization algorithms with constant step-sizes. When the global objective function is strongly convex, we prove that two variants of UGT can both achieve linear convergence. Finally, numerical experiments are taken to evaluate the performance of UGT.

\begin{appendix} \label{app}
% \section{Lemmas}
Since the proofs of \cref{the1,the2} are quite similar, we only provide the proof of the former for the sake of conciseness. We first give some lemmas which are necessary for the subsequent convergence analysis.
For CTA-UGT, define
\eqe{
    \ax^k =& \frac{1}{n}\sum_{i=1}^nx_i^k, \\
    \ag^k =& \frac{1}{n}\sum_{i=1}^ng_i^k, \\
    \bar{\nabla} f(\mx^k) =& \frac{1}{n}\sum_{i=1}^n\nabla f_i(x_i^k).
    \nonumber
}

\begin{lemma} \label{gav}
Suppose \cref{W} holds, then
\eqe{
    \ax\+ =& \ax^k - \ag^k, \\
    \ag\+ =& \alpha \bar{\nabla} f(\mx\+).
    \nonumber
}
\end{lemma}

\textit{Proof.} The first equality is obvious. Note that $\1^{\top}(\mI-\mW_2) = 0$, we have
\eqe{
    \ag\+ = \ag^k + \alpha(\bar{\nabla} f(\mx\+)-\bar{\nabla} f(\mx^k)),
    \nonumber
}
then we can obtain the second equality since $\mg^0=\alpha\nabla f(\mx^0)$.

\begin{lemma} \cite{qu2017harnessing} \label{sig}
Suppose \cref{W} holds, then
\eqe{
    \|\mW_1\mx-\1\ax\| \leq& \sigma_1\|\mx-\1\ax\|, \\
    \|\mW_2\mx-\1\ax\| \leq& \sigma_2\|\mx-\1\ax\|.
\nonumber
}
\end{lemma}

\begin{lemma} \label{ave}
\eqe{
\|\mx-\1\ax\| \leq \|\mx\|.
\nonumber
}
\end{lemma}

\textit{Proof.}
\eqe{
    \|\mx-\1\ax\|^2 &= \|\mx\|^2 + \|\1\ax\|^2 - 2\sum_{i=1}^n x_i^{\top}\ax \\
    &= \|\mx\|^2 + \|\1\ax\|^2 - 2n \ax^{\top}\ax \\
    &\leq \|\mx\|^2,
    \nonumber
}
which completes the proof.

\begin{lemma} \cite{qu2017harnessing} \label{gra}
Suppose \cref{str} holds, then
$$\|x-\alpha\nabla F(x)-x^*\| \leq \eta\|x-x^*\|,$$
where $\eta = \min\{|1-\alpha \mu|, |1-\alpha \bar{l}|\}$.
\end{lemma}

\begin{lemma} \cite{pu2021distributed} \label{rho}
Consider a nonnegative and irreducible matrix $B=[b_{ij}]\in \mathbb{R}^{3\times3}$ which satisfies that $b_{ii}<\lambda^*$, where $\lambda^*>0$ and $i=1,2,3$. Then $\rho(B) < \lambda^*$ iff $\text{det}[\lambda^*I-B]>0$.
\end{lemma}

\textit{Proof of Theorem 1.} According to CTA-UGT and \cref{gav}, we have
\eqe{ \label{bx}
    &\|\mx\+-\1\ax\+\| \\
    &= \|\mW_1\mx^k-\1\ax^k-(\mg^k-\1\ag^k)\| \\
    &\leq \sa\|\mx-\1\ax\| + \|\mg^k-\1\ag^k\|,
}
where the inequality holds due to \cref{sig}.

Notice that
\eqe{
    &\|(\mI-\mW_2)\mW_1\mx^k\| \\
    &= \|\mW_2\mW_1\mx^k-\1\ax^k -(\mW_1\mx^k- \1\ax^k)\| \\
    &\leq (1+\sigma_2)\|\mW_1\mx^k-\1\ax^k\| \\
    &\leq (1+\sigma_2)\sigma_1\|\mx^k-\1\ax^k\|,
    \nonumber
}
where the two inequalities hold because of \cref{W,sig}, then we have
\eqe{ \label{15}
    &\|\mg\+-\1\ag\+\| \\
    &\leq \|\mW_2\mg^k-\1\ag^k\| + \beta\|(\mI-\mW_2)\mW_1\mx^k\| \\
    &\quad + \|\alpha(\nabla f(\mx\+)-\nabla f(\mx^k)) - \1(\ag\+-\ag^k)\| \\
    &\leq \sigma_2\|\mg^k-\1\ag^k\| + (1+\sigma_2)\sigma_1\|\mx^k-\1\ax^k\| \\
    &\quad + \|\alpha(\nabla f(\mx\+)-\nabla f(\mx^k))\| \\
    &\leq \sigma_2\|\mg^k-\1\ag^k\| + (1+\sigma_2)\sigma_1\|\mx^k-\1\ax^k\| \\
    &\quad + \alpha l\|\mx\+-\mx^k\|,
    \nonumber
}
where the second inequality holds due to \cref{gav,ave} and the last inequality holds since $f$ is $l$-smooth.
In addition, we have
\eqe{ \label{16}
    &\|\mx\+-\mx^k\| \\
    &= \|\mW_1\mx^k-\1\ax^k-(\mx^k-\1\ax^k)-\mg^k\| \\
    &\leq (1+\sigma_1)\|\mx^k-\1\ax^k\| + \|\mg^k\|
    \nonumber
}
and
\eqe{ \label{17}
    \|\mg^k\| &\leq \|\mg^k-\1\ag^k\| + \|\1\ag^k-\1\alpha \nabla F(\ax^k)\| \\
    &\quad + \alpha \|\1\nabla F(\ax^k)-\1\nabla F(x^*)\| \\
    &\leq \|\mg^k-\1\ag^k\| + \alpha l\|\mx^k-\1\ax^k\| + \alpha l\sqrt{n}\|\ax^k-x^*\|,
    \nonumber
}
where the letter holds because
\eqe{  \label{18}
    \|\ag^k-\alpha\nabla F(\ax^k)\| &= \alpha\left\|\frac{1}{n}\sum_{i=1}^n\nabla f_i(x_i^k)-\nabla f_i(\ax^k)\right\| \\
    &\leq \alpha l\sum_{i=1}^n \frac{\|x_i^k-\ax^k\|}{n} \\
    &\leq \alpha l\sqrt{\sum_{i=1}^n \frac{\|x_i^k-\ax^k\|^2}{n}} \\
    &\leq \frac{\alpha l}{\sqrt{n}}\|\mx^k-\1\ax^k\|.
    % \nonumber
}
% Combining \cref{15}, \cref{16}, and \cref{17}, we can obtain that
It follows that
\eqe{ \label{19}
    &\|\mg\+-\1\ag\+\| \\
    &\leq [(1+\sigma_2)\sigma_1+\alpha l(1+\sigma_1+\alpha l)]\|\mx^k-\1\ax^k\| \\
    &\quad + (\sigma_2+\alpha l)\|\mg^k-\1\ag^k\| + (\alpha l)^2\sqrt{n}\|\ax^k-x^*\|.
}

Based on \cref{gav}, we have
\eqe{ \label{20}
    &\|\ax\+-x^*\| \\
    &\leq \|\ax^k-\alpha\nabla F(\ax^k)-x^*\| + \|\alpha\nabla F(\ax^k) - \ag^k\| \\
    &\leq \eta\|\ax^k-x^*\| + \frac{\alpha l}{\sqrt{n}}\|\mx^k-\1\ax^k\|,
}
where the last inequality holds because of \cref{gra,18}.

Combining \cref{bx,19,20}, we can obtain that
\eqe{
    \me\+ \leq G_1(\alpha, \beta)\me^k,
    \nonumber
}
where
$$\me^k = \left[\begin{array}{c}
\|\mx^k-\1\ax^k\| \\
\|\mg^k-\1\ag^k\| \\
\sqrt{n}\|\ax^k-x^*\|
\end{array}\right].$$

If $G_1$ is nonnegative, i.e., all entities of $G_1$ are nonnegative, we have
$$\me\+ \leq \rho(G_1)^k\me^0,$$
which means that CTA-UGT converges linearly if $\rho(G_1)<1$.

According to \cref{rho}, to guarantee $\rho(G_1)<1$, we first need to assure $0 \leq g_{ii}<1$, which holds when $0<\alpha<\frac{1-\sigma_2}{l}<\frac{1}{l}$, since $\sigma_1<1$ and $\eta = 1-\alpha\mu$ if $\alpha<\frac{1}{l}<\frac{1}{\bar{l}}$. Furthermore, we must guarantee that $\text{det}[I-G_1]>0$. Note that
\eqe{
    \text{det}[I-G_1] &= -(\alpha l)^3 + (1-\eta)[-(\alpha l)^2-2\alpha l \\
    &\quad + (1-\sigma_1)(1-\sigma_2) -\beta\sigma_1(1+\sigma_2)],
    \nonumber
}
hence it must hold that
\eqe{ \label{b12}
    (\alpha l)^3 + \alpha\mu[(\alpha l)^2+2\alpha l + \beta\sigma_1(1+\sigma_2) \\
    \quad - (1-\sigma_1)(1-\sigma_2)] < 0,
    \nonumber
}
which is equivalent to
$$\alpha^2 l^3 + \mu[(\alpha l)^2+2\alpha l + \beta\sigma_1(1+\sigma_2) - (1-\sigma_1)(1-\sigma_2)] < 0,$$
rearranging it gives that
\eqe{ \label{las}
    \alpha^2(l^3+\mu l^2) + 2\mu l\alpha + \mu[\beta\sigma_1(1+\sigma_2) \\
    \quad - (1-\sigma_1)(1-\sigma_2)] < 0.
    % \nonumber
}
Recall the expressions of $a$, $b$, and $c$ in \cref{the1} and $\beta<\frac{(1-\sigma_1)(1-\sigma_2)}{\sigma_1(1+\sigma_2)}$, then \cref{las} can be guaranteed by letting
$$\alpha <\frac{-b+\sqrt{b^2-4ac}}{2a},$$
% where $a$, $b$, and $c$ are defined in \cref{theo3}.
which finishes the proof.

% As mentioned before, the proof of \cref{the2} is almost the same with the one of \cref{the1}, which is thereby omitted.

\end{appendix}

% \begin{appendix} \label{app}
% \end{appendix}

\bibliographystyle{ieeetr}
 %\biboptions{square,numbers,sort&compress}
\bibliography{bibfile}

\begin{thebibliography}{10}

\bibitem{yang2019survey}
T.~Yang, X.~Yi, J.~Wu, Y.~Yuan, D.~Wu, Z.~Meng, Y.~Hong, H.~Wang, Z.~Lin, and
  K.~H. Johansson, ``A survey of distributed optimization,'' {\em Annual
  Reviews in Control}, vol.~47, pp.~278--305, 2019.

\bibitem{nedic2009distributed}
A.~Nedic and A.~Ozdaglar, ``Distributed subgradient methods for multi-agent
  optimization,'' {\em IEEE Transactions on Automatic Control}, vol.~54, no.~1,
  pp.~48--61, 2009.

\bibitem{shi2015extra}
W.~Shi, Q.~Ling, G.~Wu, and W.~Yin, ``Extra: An exact first-order algorithm for
  decentralized consensus optimization,'' {\em SIAM Journal on Optimization},
  vol.~25, no.~2, pp.~944--966, 2015.

\bibitem{ling2015dlm}
Q.~Ling, W.~Shi, G.~Wu, and A.~Ribeiro, ``Dlm: Decentralized linearized
  alternating direction method of multipliers,'' {\em IEEE Transactions on
  Signal Processing}, vol.~63, no.~15, pp.~4051--4064, 2015.

\bibitem{li2019decentralized}
Z.~Li, W.~Shi, and M.~Yan, ``A decentralized proximal-gradient method with
  network independent step-sizes and separated convergence rates,'' {\em IEEE
  Transactions on Signal Processing}, vol.~67, no.~17, pp.~4494--4506, 2019.

\bibitem{yuan2018exact}
K.~Yuan, B.~Ying, X.~Zhao, and A.~H. Sayed, ``Exact diffusion for distributed
  optimization and learning—part i: Algorithm development,'' {\em IEEE
  Transactions on Signal Processing}, vol.~67, no.~3, pp.~708--723, 2018.

\bibitem{nedic2017achieving}
A.~Nedic, A.~Olshevsky, and W.~Shi, ``Achieving geometric convergence for
  distributed optimization over time-varying graphs,'' {\em SIAM Journal on
  Optimization}, vol.~27, no.~4, pp.~2597--2633, 2017.

\bibitem{qu2017harnessing}
G.~Qu and N.~Li, ``Harnessing smoothness to accelerate distributed
  optimization,'' {\em IEEE Transactions on Control of Network Systems},
  vol.~5, no.~3, pp.~1245--1260, 2017.

\bibitem{di2016next}
P.~Di~Lorenzo and G.~Scutari, ``Next: In-network nonconvex optimization,'' {\em
  IEEE Transactions on Signal and Information Processing over Networks},
  vol.~2, no.~2, pp.~120--136, 2016.

\bibitem{scutari2019distributed}
G.~Scutari and Y.~Sun, ``Distributed nonconvex constrained optimization over
  time-varying digraphs,'' {\em Mathematical Programming}, vol.~176, no.~1,
  pp.~497--544, 2019.

\bibitem{xu2015augmented}
J.~Xu, S.~Zhu, Y.~C. Soh, and L.~Xie, ``Augmented distributed gradient methods
  for multi-agent optimization under uncoordinated constant stepsizes,'' in
  {\em 2015 54th IEEE Conference on Decision and Control (CDC)},
  pp.~2055--2060, IEEE, 2015.

\bibitem{nedic2017geometrically}
A.~Nedi{\'c}, A.~Olshevsky, W.~Shi, and C.~A. Uribe, ``Geometrically convergent
  distributed optimization with uncoordinated step-sizes,'' in {\em 2017
  American Control Conference (ACC)}, pp.~3950--3955, IEEE, 2017.

\bibitem{jakovetic2018unification}
D.~Jakoveti{\'c}, ``A unification and generalization of exact distributed
  first-order methods,'' {\em IEEE Transactions on Signal and Information
  Processing over Networks}, vol.~5, no.~1, pp.~31--46, 2018.

\bibitem{sundararajan2019canonical}
A.~Sundararajan, B.~Van~Scoy, and L.~Lessard, ``A canonical form for
  first-order distributed optimization algorithms,'' in {\em 2019 American
  Control Conference (ACC)}, pp.~4075--4080, IEEE, 2019.

\bibitem{alghunaim2020decentralized}
S.~A. Alghunaim, E.~Ryu, K.~Yuan, and A.~H. Sayed, ``Decentralized proximal
  gradient algorithms with linear convergence rates,'' {\em IEEE Transactions
  on Automatic Control}, 2020.

\bibitem{xu2021distributed}
J.~Xu, Y.~Tian, Y.~Sun, and G.~Scutari, ``Distributed algorithms for composite
  optimization: Unified framework and convergence analysis,'' {\em IEEE
  Transactions on Signal Processing}, 2021.

\bibitem{alghunaim2019linearly}
S.~A. Alghunaim, K.~Yuan, and A.~H. Sayed, ``A linearly convergent proximal
  gradient algorithm for decentralized optimization,'' {\em Advances In Neural
  Information Processing Systems 32 (Nips 2019)}, vol.~32, no.~CONF, 2019.

\bibitem{li2022implicit}
J.~Li and H.~Su, ``Implicit tracking-based distributed constraint-coupled
  optimization,'' {\em arXiv preprint arXiv:2201.07627}, 2022.

\bibitem{zhu2010discrete}
M.~Zhu and S.~Mart{\'\i}nez, ``Discrete-time dynamic average consensus,'' {\em
  Automatica}, vol.~46, no.~2, pp.~322--329, 2010.

\bibitem{kia2015distributed}
S.~S. Kia, J.~Cort{\'e}s, and S.~Mart{\'\i}nez, ``Distributed convex
  optimization via continuous-time coordination algorithms with discrete-time
  communication,'' {\em Automatica}, vol.~55, pp.~254--264, 2015.

\bibitem{kia2015dynamic}
S.~S. Kia, J.~Cort{\'e}s, and S.~Martinez, ``Dynamic average consensus under
  limited control authority and privacy requirements,'' {\em International
  Journal of Robust and Nonlinear Control}, vol.~25, no.~13, pp.~1941--1966,
  2015.

\bibitem{pu2021distributed}
S.~Pu and A.~Nedi{\'c}, ``Distributed stochastic gradient tracking methods,''
  {\em Mathematical Programming}, vol.~187, no.~1, pp.~409--457, 2021.

\end{thebibliography}

\end{document}